\documentclass[11pt]{report}
\usepackage{amsfonts}
 \normalbaselines

\newcommand{\cf}[1]{\mbox{\boldmath${#1}$}}
\newcommand{\E}[0]{{\rm e}}
\newcommand{\I}[0]{{\rm i}}
\newcommand{\dis}[1]{\displaystyle{#1}}
\newcommand{\text}[1]{\textstyle{#1}}
\newcommand{\scr}[1]{\scriptstyle{#1}}

\newcommand{\foot}[1]{\footnotesize{#1}}

\usepackage{wrapfig}

\usepackage[
    colorlinks=true,
    urlcolor=black,
    filecolor=black,
    linkcolor=black,
    anchorcolor=black,
    menucolor=black,
    citecolor=black,
    pdftitle={Heptagon},
    pdfauthor={Wuensche},
    pdfsubject={Approach to a Proof of the Riemann Hypothesis by the Second Mean-Value Theorem of Calculus},
    pdfkeywords={Approach to a Proof of the Riemann Hypothesis by the Second Mean-Value Theorem of Calculus},
    pdfproducer={pdfLaTeX},
    pdfpagelabels=true,
    hyperindex=true,
    plainpages=false,
    hypertexnames=false,
    linktocpage=true,
    pdfpagemode=UseThumbs, 
    pdfstartview=FitH, 
    pdfview=FitH,
    pdfmenubar=true,
    bookmarks=true,
    bookmarksnumbered=true,
    bookmarksopen=true,
    bookmarksopenlevel=1,
    bookmarksnumbered=true
]{hyperref}
\usepackage[pdftex]{graphicx}
\usepackage[pdftex]{thumbpdf}
\usepackage{xcolor}

\makeindex

\newcounter{saveeqn}
\newcommand{\aeqn}{\setcounter{saveeqn}{\value{equation}}%
\setcounter{equation}{0}%
\renewcommand{\theequation}{\mbox{A.\arabic{equation}}}}

\newcommand{\beqn}{\setcounter{saveeqn}{\value{equation}}%
\setcounter{equation}{0}%
\renewcommand{\theequation}{\mbox{B.\arabic{equation}}}}

\begin{document}
\bibliographystyle{unsrt}

\vbox {\vspace{0mm}}
\begin{center}
{\LARGE \bf  Approach to a Proof of the Riemann Hypothesis \\[4mm]by the
Second Mean-Value Theorem of Calculus}\\[4mm]
\end{center}

\begin{center}
{Alfred W\"{u}nsche}
\\[1mm]
formerly: {\it Humboldt-Universit\"{a}t Berlin,\\
Institut f\"{u}r Physik, Nichtklassische Strahlung (MPG), \\
Newtonstrasse 15, 12489 Berlin, Germany} \\
e-mail:$\;$ alfred.wuensche@physik.hu-berlin.de \\

\end{center}

{\bf{Abstract}}
{\small{
{\begin{quote}
By the second mean-value theorem of calculus (Gauss-Bonnet theorem) we prove that the class of functions ${\mit \Xi}(z)$ with an integral representation of the form $\int_{0}^{+\infty}du\,{\mit \Omega}(u)\,{\rm ch}(uz)$ with a real-valued function ${\mit \Omega}(u) \ge 0$ which is non-increasing and decreases in infinity more rapidly than any exponential functions $\exp\left(-\lambda u\right),\,\lambda >0$ possesses zeros only on the imaginary axis. The Riemann zeta function $\zeta(s)$ as it is known can be related to an entire function $\xi(s)$ with the same non-trivial zeros as $\zeta(s)$. Then after a trivial argument displacement $s\leftrightarrow z=s-\frac{1}{2}$ we relate it to a function ${\mit \Xi}(z)$ with a representation of the form ${\mit \Xi}(z)=\int_{0}^{+\infty}du\,{\mit \Omega}(u)\,{\rm ch}(uz)$ where ${\mit \Omega}(u)$ is rapidly decreasing in infinity and satisfies all requirements necessary for the given proof of the position of its zeros on the imaginary axis $z=\I y$ by the second mean-value theorem. Besides this theorem we apply the Cauchy-Riemann differential equation in an integrated operator form derived in the Appendix B. All this means that we prove a theorem for zeros of ${\mit \Xi}(z)$ on the imaginary axis $z=\I y$ for a whole class of function ${\mit \Omega}(u)$ which includes in this way the proof of the Riemann hypothesis. This whole class includes, in particular, the modified Bessel functions ${\rm I}_{\nu}(z)$ for which it is known that their zeros lie on the imaginary axis and which affirms our conclusions that we intend to publish at another place. In the same way a class of almost-periodic functions to piece-wise constant nonincreasing functions ${\rm \Omega}(u)$ belong also to this case.
At the end we give shortly an equivalent way of a more formal description of the obtained results using the Mellin transform of functions with its variable substituted by an operator.

\end{quote}
}}

\setcounter{chapter}{1}
\section*{1. Introduction}

The Riemann zeta function $\zeta(s)$ which basically was known already to Euler establishes the most important link between number theory and analysis. The proof of the Riemann hypothesis is a longstanding problem since it was formulated by Riemann \cite{rie} in 1859. The Riemann hypothesis is the conjecture that all nontrivial zeros of the Riemann zeta function $\zeta(s)$ for complex $s=\sigma+\I t$ are positioned on the line $s=\frac{1}{2}+\I t$ that means on the line parallel to the imaginary axis through real value $\sigma = \frac{1}{2}$ in the complex plane and in extension that all zeros are simple zeros \cite{whitt,titch,chand,edw,parschin,patt,cart,ivic,havil,riben,borw,apost0,weisst,
lang,
spekt,stew} (with extensive lists of references in some of the cited sources, e.g., \cite{chand,edw,ivic,borw,weisst}). The book of
Edwards \cite{edw} is one of the best older sources concerning most problems connected with the Riemann zeta function. There are also mathematical tables and chapters in works about Special functions which contain information about the Riemann zeta function and about number analysis,
e.g., Whittaker and Watson \cite{whitt} (chap. 13), Bateman and Erd\'{e}lyi \cite{bate1} (chap. 1) about zeta functions, \cite{bate3} (chap. 17) about number analysis, and Apostol \cite{apost,nist} (chaps. 25 and 27). The book of Borwein, Choi, Rooney and Weirathmueller \cite{borw} gives on the first 90 pages a short account about  achievements concerning the Riemann hypothesis and its consequences for number theory and on the following about 400 pages it reprints important original papers and expert witnesses in the field.
Riemann has put aside the search for a proof of his hypothesis 'after some fleeting vain attempts' and emphasizes that 'it is not necessary for the immediate objections of his investigations' \cite{rie} (see \cite{edw}).
The Riemann hypothesis was taken by Hilbert as the 8-th problem in his representation of 23 fundamental unsolved problems in pure mathematics and axiomatic physics in a lecture hold on 8 August in 1900 at the Second Congress of Mathematicians in Paris \cite{hilb,linnik}.
The vast experience with the Riemann zeta function in the past and the progress in numerical calculations of the zeros (see, e.g., \cite{edw,havil,riben,spekt,stew,conrey,bomb}) which all confirmed the Riemann hypothesis suggest that it should be true corresponding to the opinion of most of the specialists in this field but not of all specialists (arguments for doubt are discussed in \cite{ivic1}).

The Riemann hypothesis is very important for prime number theory and a number of consequences is derived under the unproven assumption that it is true. As already said a main role plays a function $\zeta(s)$ which was known already to Euler for real variables $s$ in its product representation (Euler product) and in its series representation and was continued to the whole complex $s$-plane by Riemann and is now called Riemann zeta function. The Riemann hypothesis as said is the conjecture that all nontrivial zeros of the zeta function $\zeta(s)$ lie on the axis parallel to the imaginary axis and intersecting the real axis at $s=\frac{1}{2}$. For the true hypothesis the representation of the Riemann zeta function after exclusion of its only singularity at $s=1$ and of the trivial zeros at $s=-2n,\,(n=1,2,\ldots)$ on the negative real axis is possible by a Weierstrass product with factors which only vanish on the critical line $\sigma=\frac{1}{2}$. The function which is best suited for this purpose is the so-called xi function $\xi(s)$ which is closely related to the zeta function $\zeta(s)$ and which was also introduced by Riemann \cite{rie}. It contains all information about the nontrivial zeros and only the exact positions of the zeros on this line are not yet given then by a closed formula which, likely, is hardly to find explicitly but an approximation for its density was conjectured already by Riemann \cite{rie} and proved by von Mangoldt \cite{mang2}.
The "(pseudo)-random" character of this distribution of zeros on the critical line remembers somehow the "(pseudo)-random" character of the distribution of primes where one of the differences is that the distribution of primes within the natural numbers becomes less dense with increasing integers whereas the distributions of zeros of the zeta function on the critical line becomes more dense with higher absolute values with slow increase and approaches to a logarithmic function in infinity.

There are new ideas for analogies to and application of the Riemann zeta function in other regions of mathematics and physics. One direction is the theory of random matrices \cite{spekt,conrey} which shows analogies in their eigenvalues to the distribution of the nontrivial zeros of the Riemann zeta function. Another interesting idea founded by Voronin \cite{voronin} (see also \cite{spekt,steud}) is the universality of this function in the sense that each holomorphic function without zeros and poles in a certain circle with radius less $\frac{1}{2}$ can be approximated with arbitrary required accurateness in a small domains of the zeta function to the right of the critical line within $\frac{1}{2} \le s \le 1$. An interesting idea is elaborated in articles of Neuberger, Feiler, Maier and Schleich \cite{schl1,schl2}. They consider a simple first-order ordinary differential equation with a real variable $t$ (say the time) for given arbitrary analytic functions $f(z)$ where the time evolution of the function for every point $z$ finally transforms the function in one of the zeros $f(z)=0$ of this function in the complex $z$-plane and illustrate this process graphically by flow curves which they call Newton flow and which show in addition to the zeros the separatrices of the regions of attraction to the zeros. Among many other functions they apply this to the Riemann zeta function $\zeta(z)$ in different domains of the complex plane. Whether, however, this may lead also to a proof of the Riemann hypothesis is more than questionable.

Number analysis defines some functions of a continuous variable, for example, the number of primes $\pi(x)$ less a given real number $x$ which last is connected with the discrete prime number distribution (e.g., \cite{titch,chand,edw,patt,ivic,riben}) and establishes the connection to the Riemann zeta function $\zeta(s)$. Apart from the product representation of the Riemann zeta function the representation by
a type of series which is now called Dirichlet series was already known to Euler. With these Dirichlet series in number theory are connected some discrete functions over the positive integers $n=1,2,\ldots$ which play a role as coefficients in these series and are called arithmetic functions (see, e.g., Chandrasekharan \cite{chand} and Apostol \cite{apost0}). Such functions are the M\"{o}bius function $\mu(n)$ and the Mangoldt function ${\rm \Lambda}(n)$ as the best known ones. A short representation of the connection of the Riemann zeta function to number analysis and of some of the functions defined there became now standard in many monographs about complex analysis (e.g., \cite{lang}).

Our means for the proof of the Riemann hypothesis in present article are more conventional and "old-fashioned" ones, i.e. the Real Analysis and the Theory of Complex Functions which were developed already for a long time.
The most promising way for a proof of the Riemann hypothesis as it seemed to us in past is via the already mentioned entire function $\xi(s)$ which is closely related to the Riemann zeta function $\zeta(s)$. It contains all important elements and information of the last but excludes its trivial zeros and its only singularity and, moreover, possesses remarkable symmetries which facilitate the work with it compared with the Riemann zeta function. This function $\xi(s)$ was already introduced by Riemann \cite{rie} and dealt with, for example, in the classical books of Titchmarsh \cite{titch}, Edwards \cite{edw} and in almost all of the sources cited at the beginning. Present article is mainly concerned with this xi function $\xi(s)$ and its investigation in which, for convenience, we displace the imaginary axis by $\frac{1}{2}$ to the right that means to the critical line and call this Xi function ${\mit\Xi}(z)$ with $z=x+\I y$. We derive some representations for it among them novel ones and discuss its properties, including its derivatives, its specialization to the critical line and some other features. We make an approach to this function via the second mean value theorem of analysis (Gauss-Bonnet theorem, e.g., \cite{courant,widder}) and then we apply an operator identity for analytic functions which is derived in Appendix B and which is equivalent to a somehow integrated form of the Cauchy-Riemann equations. This among other not so successful trials (e.g., via moments of function ${\mit \Omega}(u)$) led us finally to a proof of the Riemann hypothesis embedded into a proof for a more general class of functions.

Our approach to a proof of the Riemann hypothesis in article in rough steps is as follows:\\
First we shortly represent the transition from the Riemann zeta function $\zeta(s)$ of complex variable $s=\sigma +\I t$ to the xi function $\xi(s)$ introduced already by Riemann and derive for it by means of the Poisson summation formula a representation which is convergent in the whole complex plane (Section 2 with main formal part in Appendix A).
Then we displace the imaginary axis of variable $s$ to the critical line at $s=\frac{1}{2}+\I t$ by $s\rightarrow z = s-\frac{1}{2}$ that is purely for convenience of further working with the formulae. However, this has also the desired subsidiary effect that it brings us into the fairway of the complex analysis usually represented with the complex variable $z=x+\I y$. The transformed $\xi(s)$ function is called ${\mit \Xi}(z)$ function.

The function ${\mit \Xi}(z)$ is represented as an integral transform of a real-valued function ${\mit \Omega}(u)$ of the real variable $u$ in the form ${\mit \Xi}(z)= \int_{0}^{+\infty}du\,{\mit \Omega}(u)\,{\rm ch}(u z)$ which is related to a Fourier transform (more exactly to Cosine Fourier transform). If the Riemann hypothesis is true then we have to prove that all zeros of the function ${\mit \Xi}(z)$ occur for $x=0$.

To the Xi function in mentioned integral transform we apply the second mean-value theorem of real analysis first on the imaginary axes and discuss then its extension from the imaginary axis to the whole complex plane. For this purpose we derive in Appendix B in operator form relations of the complex mean-value parameter in our application of the second mean-value theorem to this parameter on the imaginary axis which are equivalents in integral form to the Cauchy-Riemann equations in differential form and apply this in specific form to the Xi function (Sections 3 and 4).

Then in Section 5 we accomplish the proof with the discussion and solution of the two most important equations (\ref{zeroscond1}) and (\ref{zeroscond2}) for the last and decisive stage of the proof. These two equations are derived before in preparation of this last stage of the proof. From these equations it is seen that the obtained two real equations admit zeros of the Xi function only on the imaginary axis. This proves the Riemann hypothesis by the equivalence of the Riemann zeta function $\zeta(s)$ to the Xi function ${\mit \Xi}(z)$ and embeds it into a whole class of functions with similar properties and positions of their zeros.

The Sections 6-7 serve for illustrations and graphical representations of the specific parameters (e.g., mean-value parameters) for the Xi function to the Riemann hypothesis and for other functions which in our proof by the second mean-value problem are included for the existence of zeros only on the imaginary axis. This is in particular the whole class of modified Bessel functions ${\rm I}_{\nu}(z), \; \left(-\frac{1}{2} < \nu < +\infty\right)$ with real indices $\nu$ which possess zeros only on the imaginary axis $y$ and where a proof by means of the differential equations exists.

\setcounter{chapter}{2}
\setcounter{equation}{0}
\section*{2. From Riemann zeta function ${\cf{\zeta(s)}}$ to related xi function ${\cf{\xi(s)}}$ and its argument displacement to function ${\cf{{\mit \Xi(z)}}} $}

In this Section we represent the known transition from the Riemann zeta function $\zeta(s)$ to a function $\xi(s)$ and finally to a function ${\mit \Xi}(z)$ with displaced complex variable $s\rightarrow z=s-\frac{1}{2}$ for rational effective work and establish some of the basic representations of these functions, in particular, a kind of modified Cosine Fourier transformations of a function ${\mit \Omega}(u)$ to the function ${\mit \Xi}(z)$.

As already expressed in the Introduction, the most promising way for a proof of the Riemann hypothesis as it seems to us is the way via a certain integral representation of the related xi function $\xi(s)$. We sketch here the transition from the Riemann zeta function $\zeta(s)$ to the related xi function $\xi(s)$ in a short way because, in principle, it is known and we delegate some aspects of the derivations to Appendix A

Usually, the starting point for the introduction of the Riemann zeta function $\zeta(s)$ is the following relation between the Euler product and an infinite series continued to the whole complex $s$-plane
\begin{eqnarray}
\zeta(s) \equiv \prod_{n=1}^\infty
\left(1-\frac{1}{p_n^{\,s}}\right)^{-1} = \sum_{n=1}^\infty
\frac{1}{n^s},\quad \left(\sigma \equiv {\rm Re}(s)
>1\right),\label{eulerprod}
\end{eqnarray}
where $p_n$ denotes the ordered sequence of primes ($p_1=2,p_2=3,p_3=5,\ldots$). The transition from the product formula to the sum representation in (\ref{eulerprod}) via transition to the Logarithm of $\zeta(s)$ and Taylor series expansion of the factors $\log\left(1-\frac{1}{p_n^{\,s}}\right)^{-1}$ in powers of $\frac{1}{p_n^{\,s}}$ using the uniqueness of the prime-number decomposition is well known and due to Euler in 1737. It leads to a special case of a kind of series later introduced and investigated in more general form and called Dirichlet series. The Riemann zeta function $\zeta(s)$ can be analytically continued into the whole complex plane to a meromorphic function that was made and used by Riemann. The sum in (\ref{eulerprod}) converges uniformly for complex variable
$s=\sigma +\I t$ in the open semi-planes with arbitrary $\sigma > 1$ and arbitrary $t$. The only singularity of the function $\zeta(s)$ is a simple pole at $s=1$ with residue $1$ that we discuss below.

The product form (\ref{eulerprod}) of the zeta function $\zeta(s)$ shows that it involves all prime numbers $p_n$ exactly one times and therefore it contains information about them in a coded form. It proves to be possible to regain information about the prime number distribution from this function. For many purposes it is easier to work with meromorphic and, moreover, entire functions than with infinite sequences of numbers but in first case one has to know the properties of these functions which are determined by their zeros and their singularities together with their multiplicity.

From the well-known integral representation of the Gamma function
\begin{eqnarray}
{\rm \Gamma}(z) &=& \int_{0}^{+\infty}dt\,t^{z-1}\E^{-t},\quad
\left({\rm Re}(z) > 0\right),\label{gammadef}
\end{eqnarray}
follows by the substitutions $t=n^{\mu}x,\,\mu z =s$ with an appropriately fixed parameter $\mu > 0$ for arbitrary natural numbers $n$
\begin{eqnarray}
\frac{1}{n^s} &=& \frac{1}{{\rm
\Gamma}\left(\frac{s}{\mu}\right)}\int_{0}^{+\infty}dx\,x^{\frac{s}{\mu}-1}
\E^{-n^{\mu}x},\quad \left({\rm
Re}\left(\frac{s}{\mu}\right)>0\right).\label{gammadefcons}
\end{eqnarray}
Inserting this into the sum representation (\ref{eulerprod}) and changing the order of summation and integration, we obtain for choice $\mu=1$ of the parameter using the sum evaluation of the geometric series
\begin{eqnarray}
\zeta(s) &=& \frac{1}{{\rm
\Gamma}(s)}\int_{0}^{+\infty}dp\,\frac{p^{s-1}}{\E^{p}-1},\quad
\left({\rm Re}(s) > 1\right),\label{zetarep1}
\end{eqnarray}
and for choice $\mu=2$ with substitution $p=\pi q^2$ of the integration variable (see \cite{rie} and, e.g.,
\cite{titch,chand,edw,patt,ivic})
\begin{eqnarray}
\zeta(s) &=& \frac{\pi^{\frac{s}{2}}s}{{\rm
\Gamma}\left(\frac{s}{2}+1\right)}\int_{0}^{+\infty}dq\,q^{s-1}
\sum_{n=1}^\infty\exp\left(-\pi n^2q^2\right),\quad \left({\rm
Re}(s) > 1\right). \label{zetarep2}
\end{eqnarray}
Other choice of $\mu$ seems to be of lesser importance. Both representations (\ref{zetarep1}) and (\ref{zetarep2}) are closely related to a Mellin transform $\hat{f}(s)$ of a function $f(t)$ which together with its inversion is generally defined by (e.g., \cite{lang,erd,zem,brpr,bert})
\begin{eqnarray}
f(t) \rightarrow \hat{f}(s) &\equiv & \int_{0}^{+\infty} dt\,t^{s-1}\,f(t),\quad
\Rightarrow \quad f(\lambda t) \rightarrow \frac{1}{\lambda^s}\hat{f}(s), \quad \left(\,\lambda >0\,\right),\nonumber\\ \hat{f}(s) \rightarrow f(t) &=& \frac{1}{\I 2\pi}\int_{c-\I
\infty}^{c+\I \infty}ds\,t^{-s}\,\hat{f}(s), \quad
\Rightarrow \quad \hat{f}(s-s_0) \rightarrow t^{s_0}f(t),\label{mellin}
\end{eqnarray}
where $c$ is an arbitrary real value within the convergence strip of $\hat{f}(s)$ in complex $s$-plane.
The Mellin transform $\hat{f}(s)$ of a function $f(t)$ is closely related to the Fourier transform $\tilde{\varphi}(y)$ of the function $\varphi(x) \equiv f(\E^x)$ by variable substitution $t=\E^{x}$ and $y=\I s$. Thus the Riemann zeta function $\zeta(s)$ can be represented, substantially (i.e., up to factors depending on $s$), as the Mellin transforms of the functions $f(t)=\sum_{n=1}^\infty \E^{-nt}=\frac{1}{\E^{t}-1}$ or of $f(t)=\sum_{n=1}^\infty\exp\left(-\pi n^2 t^2\right)$, respectively. The kernels of the Mellin transform are the eigenfunctions of the differential operator $t\frac{\partial}{\partial t}$ to eigenvalue $s-1$ or, correspondingly, of the integral operator $\exp\left(\alpha\, t\frac{\partial}{\partial t}\right)$ of the multiplication of the argument of a function by a factor $\E^{\alpha}$ (scaling of argument).
Both representations (\ref{zetarep1}) and (\ref{zetarep2}) can be used for the derivation of further representations of the Riemann zeta function and for the analytic continuation. The analytic continuation of the Riemann zeta function can also be obtained using the Euler-Maclaurin
summation formula for the series in (\ref{eulerprod}) (e.g., \cite{edw,riben,lang}).

Using the Poisson summation formula, one can transform the representation (\ref{zetarep2}) of the Riemann zeta function to the following form
\begin{eqnarray}
\zeta(s) &=&
\frac{\pi^{\frac{s}{2}}}{\left(\frac{s}{2}\right)!(s-1)}\left\{\frac{1}{2}-s(1-s) \int_{1}^{+\infty}dq\,\frac{q^s+q^{1-s}}{q}\sum_{n=1}^\infty \exp\left(-\pi n^2
q^2\right)\right\}.\label{zetarep3}
\end{eqnarray}
This is known \cite{rie,titch,chand,edw,patt,ivic} but for convenience and due to the importance of this representation for our purpose we give a derivation in Appendix A. From (\ref{zetarep3})
which is now already true for arbitrary complex $s$ and, therefore, is an analytic continuation of the representations (\ref{eulerprod}) or (\ref{zetarep2}) we see that the Riemann zeta function satisfies a functional equation for the transformation of the argument $s\rightarrow 1-s$. In simplest form it appears by 'renormalizing' this function via introduction of the xi function $\xi(s)$ defined by Riemann according to \cite{rie} and to \cite{edw,apost}\footnote{Riemann \cite{rie} defines it more specially for argument $s=\frac{1}{2}+\I t$ and writes it $\xi(t)$ with real $t$ corresponding to our $\xi\left(\frac{1}{2}+\I t\right)$. Our definition agrees, e.g., with Eq. (1) in Section 1.8 on p. 16 of Edwards \cite{edw} and with \cite{apost} and many others.}
\begin{eqnarray}
\xi(s) & \equiv & \frac{(s-1)\left(\frac{s}{2}\right)!} {\pi^{\frac{s}{2}}}\zeta(s),\label{xirep1}
\end{eqnarray}
and we obtain for it the following representation converging in the whole complex plane of $s$ (e.g., \cite{rie,chand,edw,ivic,patt})
\begin{eqnarray}
\xi(s) &=& \frac{1}{2}-s(1-s)\int_{1}^{+\infty}dp\,\frac{p^s+p^{1-s}}{p}\sum_{n=1}^\infty
\exp\left(-\pi n^2 p^2\right),\label{xirep2}
\end{eqnarray}
with the 'normalization'
\begin{eqnarray}
\xi(0) \;=\; \xi(1) \;=\; -\zeta(0)\;=\; \frac{1}{2}.\label{normxi}
\end{eqnarray}
For $s=\frac{1}{2}$ the xi function and the zeta function possess the (likely transcendental) values
\begin{eqnarray}
\xi\left(\frac{1}{2}\right) \;=\; - \frac{\left(\frac{1}{4}\right)!}{2\pi^{\frac{1}{4}}}\,\zeta
\left(\frac{1}{2}\right)\;=\; 0.4971207782,\quad \zeta\left(\frac{1}{2}\right) \;=\; -1.4603545088 \label{xizero}
\end{eqnarray}
Contrary to the Riemann zeta function $\zeta(s)$ the function $\xi(s)$ is an entire function. The only singularity of $\zeta(s)$ which is the simple pole at $s=1$, is removed by multiplication of $\zeta(s)$ with $s-1$ in the definition (\ref{xirep1}) and the trivial zeros of $\zeta(s)$ at $s=-2n,\,(n=1,2,\ldots)$ are also removed by its multiplication with $\left(\frac{s}{2}\right)! \equiv \,{\rm \Gamma}\left(\frac{s}{2}+1\right)$ which possesses simple poles there.

The functional equation
\begin{eqnarray}
\xi(s) &=& \xi(1-s),\label{funceqxi}
\end{eqnarray}
from which follows for the $n$-th derivatives
\begin{eqnarray}
\xi^{(n)}(s) &=& (-1)^n \xi^{(n)}(1-s),\quad \Rightarrow \quad \xi^{(2m+1)}\left(\frac{1}{2}\right) \;=\;0,\quad (n,m=0,1,2,\ldots),\label{funceqxider}
\end{eqnarray}
and which expresses that $\xi(s)$ is a symmetric function with respect to $s=\frac{1}{2}$ as it
is immediately seen from (\ref{xirep2}) and as it was first derived by Riemann \cite{rie}. It can be easily converted into the following functional equation for the Riemann zeta function $\zeta(s)$\footnote{According to Havil \cite{havil}, (p. 193), already Euler correctly conjectured this relation for the zeta function $\zeta(s)$ which is equivalent to relation (\ref{funceqxi}) for the function $\xi(s)$ but could not prove it. Only Riemann proved it first.}
\begin{eqnarray}
\zeta(s) &=& \frac{\left(2\pi\right)^{s}}{2\,(s-1)!\cos\left(\pi\frac{s}{2}\right)}\,\zeta(1-s).\label{funceqzeta}
\end{eqnarray}
Together with $\xi\left(s\right) =6 \left(\xi\left(s^*\right)\right)^*$ we find by combination with (\ref{funceqxi})
\begin{eqnarray}
\xi\left(s\right) \;=\; \xi(1-s) \;=\; \left(\xi\left(1-s^*\right)\right)^* \;=\;\left(\xi\left(s^*\right)\right)^*,\label{xisymmrel}
\end{eqnarray}
that combine in simple way, function values for 4 points $(s,1-s,1-s^*,s^*)$ of the complex plane.
Relation (\ref{xisymmrel}) means that in contrast to the function $\zeta(s)$ which is only real-valued on the real axis the function $\xi(s)$ becomes real-valued on the real axis ($s=s^*$) and on the imaginary axis ($s=-s^*$).

As a consequence of absent zeros of the Riemann zeta function $\zeta(\sigma +\I t)$ for $\sigma \equiv {\rm Re}(s) > 1$ together with the functional relation (\ref{funceqzeta}) follows that all nontrivial zeros of this function have to be within the strip $0\le \sigma \le 1$ and the Riemann hypothesis asserts that all zeros of the related xi function $\xi(s)$ are positioned on the so-called critical line $s=\frac{1}{2}+\I t,\,(-\infty < t < +\infty)$. This is, in principle, well known.

We use the functional equation (\ref{funceqxi}) for a simplification of the notations in the following considerations and displace the imaginary axis of the complex variable $s=\sigma +\I t$ from $\sigma=0$ to the value $\sigma=\frac{1}{2}$ by introducing the entire function ${\mit \Xi}(z)$ of the complex variable $z=x+\I y$ as follows
\begin{eqnarray}
{\mit \Xi}\left(z\right) &\equiv & \xi\left(\frac{1}{2}+z\right),\quad z=x+\I y = \sigma-\frac{1}{2}+\I t = s -\frac{1}{2},\label{capxidef}
\end{eqnarray}
with the 'normalization' (see (\ref{normxi}) and (\ref{xizero}))
\begin{eqnarray}
{\mit \Xi}\left(\pm\frac{1}{2}\right) \;=\; \frac{1}{2},\quad  {\mit \Xi}\left(0\right) \;=\; \xi\left(\frac{1}{2}\right) \;\approx\;  0.4971207782.\label{normcapxi}
\end{eqnarray}
following from (\ref{normxi}). Thus the full relation of the Xi function $\Xi(z)$ to the Riemann zeta function $\zeta(s)$ using definition (\ref{xirep1}) is
\begin{eqnarray}
{\mit \Xi}(z) &=& \frac{\left(z-\frac{1}{2}\right)\left(\frac{1+2z}{4}\right)!}{\pi^\frac{1+2z}{4}}
\zeta\left(\frac{1}{2}+z\right).\label{xizeta}
\end{eqnarray}
We emphasize again that the argument displacement (\ref{capxidef}) is made in the following only for convenience of notations and not for some more principal reason.

The functional equation (\ref{funceqxi}) together with (\ref{funceqxider}) becomes
\begin{eqnarray}
{\mit \Xi}\left(z\right) &=& {\mit \Xi}\left(-z\right),\quad {\mit \Xi}^{(n)}\left(z\right) \;=\; (-1)^n {\mit \Xi}^{(n)}\left(-z\right), \label{funceqcapxi}
\end{eqnarray}
and taken together with the symmetry for the transition to complex conjugated variable
\begin{eqnarray}
{\mit \Xi}\left(z\right) \;=\; {\mit \Xi}\left(-z\right) \;=\; \left({\mit \Xi}\left(-z^*\right)\right)^* \;=\; \left({\mit \Xi}\left(z^*\right)\right)^*.\label{xisymm}
\end{eqnarray}
This means that the Xi function ${\mit \Xi}\left(z\right)$ becomes real-valued on the imaginary axis $z=\I\,y$ which becomes the critical line in the new variable $z$
\begin{eqnarray}
{\mit \Xi}\left(\I y\right) \;=\; {\mit \Xi}\left(-\I y\right) \;=\; \left({\mit \Xi}\left(\I y\right)\right)^*  \;=\; \left({\mit \Xi}\left(-\I y\right)\right)^*.
\end{eqnarray}
Furthermore, the function ${\mit \Xi}(z)$ becomes a symmetrical function and a real-valued one on the real axis $z=x$
\begin{eqnarray}
{\mit \Xi}\left(x\right) = {\mit \Xi}\left(-x\right) = \left({\mit \Xi}\left(-x\right)\right)^* = \left({\mit \Xi}\left(x\right)\right)^*.
\end{eqnarray}
In contrast to this the Riemann zeta function $\zeta(s)$ the function is not a real-valued function on the
critical line $s=\frac{1}{2}+\I t$ and is real-valued but not symmetric on the real axis.
This is represented in Fig. 2.1. (calculated with "Mathematica 6" such as the further figures too).
\begin{figure}[h]
\includegraphics[width=16.0cm]{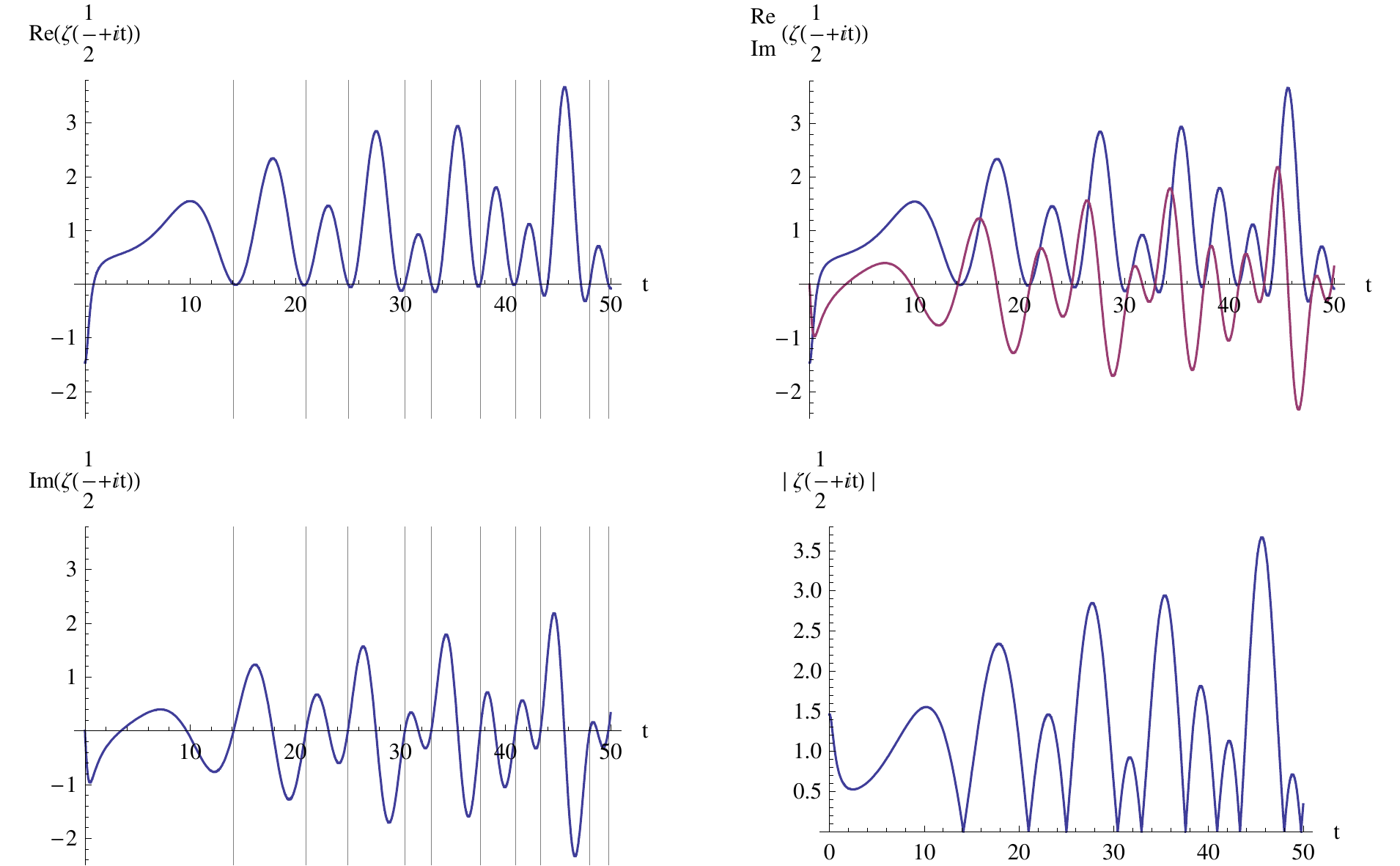}\hspace{1cm}
\caption{\foot Real and imaginary part and absolute value of Riemann zeta function on critical line. \newline \scriptsize{The position of the zeros of the whole function $\zeta\left(\frac{1}{2}+\I t\right)$ on the critical line are shown by grid lines. One can see that not all zeros of the real part are also zeros of the imaginary part and vice versa. The figures are easily to generate by program "Mathematica" and are published in similar forms already in literature.}}
\end{figure}
We see that not all of the zeros of the real part ${\rm Re}\left(\zeta\left(\frac{1}{2}+\I t\right)\right)$ are also zeros of the imaginary part ${\rm Im}\left(\zeta\left(\frac{1}{2}+\I t\right)\right)$ and, vice versa, that not all of the zeros of the imaginary part are also zeros of the real part and thus genuine zeros of the function $\zeta\left(\frac{1}{2}+\I t\right)$ which are signified by grid lines. Between two zeros of the real part which are genuine zeros of $\zeta\left(\frac{1}{2}+\I t\right)$ lies in each case (exception first interval) an additional zero of the imaginary part, which almost coincides with a maximum of the real part.

Using (\ref{xirep2}) and definition (\ref{capxidef}) we find the following representation of ${\mit \Xi(z)}$
\begin{eqnarray}
{\mit \Xi}\left(z\right) &=& \frac{1}{2}-\left(\frac{1}{4}-z^2\right) \int_{1}^{+\infty}dq\,\frac{q^z+q^{-z}}{\sqrt{q}} \sum_{n=1}^\infty \exp\left(-\pi n^2 q^2\right).\label{capxirep1}
\end{eqnarray}
With the substitution of the integration variable $q= \E^u $ (see also (\ref{yex}) in Appendix A) representation (\ref{capxirep1}) is transformed to
\begin{eqnarray}
{\mit \Xi}\left(z\right) &\equiv & \frac{1}{2}-2\left(\frac{1}{4}-z^2\right)\int_{0}^{+\infty}du\,{\rm
ch}\left(uz\right)\,\E^{\frac{u}{2}} \sum_{n=1}^\infty\exp\left(-\pi n^2 \E^{2u}\right). \label{capxirep2}
\end{eqnarray}
In Appendix A we show that (\ref{capxirep2}) can be represented as follows (see also Eq. (2) on p. 17 in \cite{edw} which possesses a similar principal form)
\begin{eqnarray}
{\mit \Xi}\left(z\right) &=& \int_{0}^{+\infty}du\,{\mit
\Omega}\left(u\right)\,{\rm ch}\left(uz\right),\label{capxirep3}
\end{eqnarray}
with the following explicit form of the function ${\mit \Omega}\left(u\right)$ of the real variable $u$
\begin{eqnarray}
{\mit \Omega}\left(u\right) &\equiv & 4\E^{\frac{u}{2}}\sum_{n=1}^\infty \pi n^2\E^{2u}\left(2\pi
n^2\E^{2u}-3\right)\exp\left(-\pi n^2 \E^{2u}\right)\;>\;0,\quad (-\infty < u < +\infty). \label{defcapomega}
\end{eqnarray}
The function ${\mit \Omega}\left(u\right)$ is symmetric
\begin{eqnarray}
{\mit \Omega}\left(u\right) = +{\mit \Omega}\left(-u\right) = {\mit
\Omega}\left(\left|u\right|\right),\label{omegasym}
\end{eqnarray}
that means it is an even function although this is not immediately seen from representation (\ref{defcapomega})\footnote{It was for us for the first time and was very surprising to meet a function where its symmetry was not easily seen from its explicit representation. However, if we substitute in (\ref{defcapomega}) $u\rightarrow -u$ and calculate and plot the part of ${\mit \Omega}(u)$ for $u\ge 0$ with the obtained formula then we need much more sum terms for the same accurateness than in case of calculation with (\ref{defcapomega}).)}. We prove this in Appendix B. Due to this symmetry, formula (\ref{capxirep3}) can be also represented by
\begin{eqnarray}
{\mit \Xi}\left(z\right) &=& \frac{1}{2}\int_{-\infty}^{+\infty}du\,{\mit
\Omega}\left(u\right)\,{\rm ch}\left(uz\right) \;=\; \frac{1}{2}\int_{-\infty}^{+\infty}du\,{\mit
\Omega}\left(u\right)\,\E^{uz}.\label{capxirep4}
\end{eqnarray}
In the formulation of the right-hand side the function ${\mit \Xi}\left(z\right)$ appears as analytic continuation of the Fourier transform of the function ${\mit \Omega}\left(u\right)$ written with imaginary argument $z=\I y$ or, more generally, with substitution $z\rightarrow \I z'$ and complex $z'$. From this follows as inversion of the integral transformation (\ref{capxirep4}) using (\ref{omegasym})
\begin{eqnarray}
{\mit \Omega}\left(u\right) &=& \frac{1}{\pi}\int_{-\infty}^{+\infty}dy\,{\mit \Xi}\left(\I y\right)\E^{-\I u y} \;=\; \frac{1}{\pi}\int_{-\infty}^{+\infty}dy\,{\mit \Xi}\left(\I y\right)\cos\left(u y\right),
\end{eqnarray}
or due to symmetry of the integrand in analogy to (\ref{capxirep3})
\begin{eqnarray}
{\mit \Omega}\left(u\right) &=& \frac{2}{\pi}\int_{0}^{+\infty}dy\,{\mit \Xi}\left(\I y\right)\cos\left(u y\right),
\end{eqnarray}
where ${\mit \Xi}\left(\I y\right)$ is a real-valued function of the variable $y$ on the imaginary axis
\begin{eqnarray}
{\mit \Xi}\left(\I y\right) &=& \int_{0}^{+\infty}du\,{\mit
\Omega}\left(u\right)\cos\left(uy\right),\quad \left({\mit \Xi}\left(\I y\right)\right)^*={\mit \Xi}\left(\I y\right),\label{xisym}
\end{eqnarray}
due to (\ref{capxirep3}).

A graphical representation of the function $\Omega(u)$ and of its first derivatives $\Omega^{(1)}(u),(n=1,2,3)$ is given in Fig. 2.2.
\begin{figure}[h]
\includegraphics[width=16.0cm]{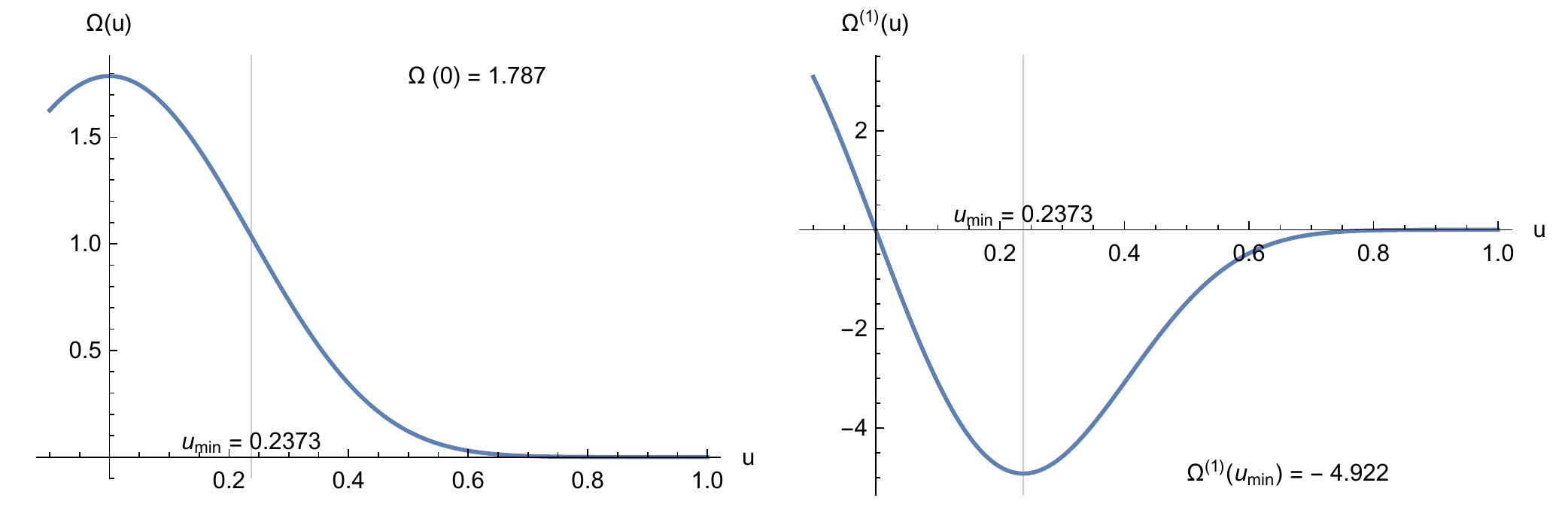}\hspace{1cm}
\caption{\foot Function $\Omega(u)$ and its first derivative $\Omega^{(1)}(u),\,(n=1,2,3)$ (see (\ref{capxirep3}) and (\ref{capxirep5})).\newline \scriptsize{The function ${\mit \Omega}(u)$ is positive for $ 0 \le u < +\infty $ and since its first derivative ${\mit \Omega}^{(1)}(u)$ is negative for $0 < u < +\infty$ the function ${\mit \Omega}(u)$ is monotonically decreasing on the real positive axis. It vanishes in infinity more rapidly than any exponential function with a polynomial in the exponent.}}
\end{figure}
The function ${\mit \Omega}\left(u\right)$ is monotonically decreasing for $0\le u < +\infty $ due to the non-positivity of its first derivative ${\mit \Omega}^{(1)}(u) \equiv \frac{\partial {\mit \Omega(u)}}{\partial u}$ which explicitly is (see also Appendix A)
\begin{eqnarray}
{\mit \Omega}^{(1)}(u) &=& -2 \E^{\frac{u}{2}}\sum_{n=1}^\infty \pi n^2\E^{2u}\left(8\left(\pi n^2\E^{2u}\right)^2-30\pi n^2\E^{2u}+15\right)\exp\left(-\pi n^2 \E^{2u}\right) \nonumber\\ &\le & 0,\qquad \left(0\le u < +\infty\right),\label{dercapomega}
\end{eqnarray}
with one relative minimum at $u_{\rm min}=0.237266$ of depth ${\mit \Omega}^{(1)}\left(u_{\rm min}\right)= -4.92176$. Moreover, it is very important for the following that due to presence of factors $\exp\left(-\pi n^2 \E^{2u}\right)$ in the sum terms in (\ref{defcapomega}) or in (\ref{dercapomega}) the functions ${\mit \Omega}(u)$ and ${\mit \Omega}^{(1)}(u)$ and all their higher derivatives are very rapidly decreasing for $u\rightarrow +\infty$, more rapidly than any exponential function with a polynomial of $u$ in the argument. In this sense the function ${\mit \Omega}(u)$ is more comparable with functions of finite support which vanish from a certain $u \ge u_0$ on than with any exponentially decreasing function.
From (\ref{omegasym}) follows immediately that the function ${\mit \Omega}^{(1)}(u)$ is antisymmetric
\begin{eqnarray}
{\mit \Omega}^{(1)}\left(u\right) = -{\mit \Omega}^{(1)}\left(-u\right) = \frac{u}{|u|}\frac{\partial}{\partial |u|}{\mit \Omega}\left(\left|u\right|\right),\quad \Rightarrow \quad {\mit \Omega}^{(1)}\left(0\right) =0,\label{omega1sym}
\end{eqnarray}
that means it is an odd function.

It is known that smoothness and rapidness of decreasing in infinity of a function change their role in Fourier transformations. As the Fourier transform of the smooth (infinitely continuously differentiable) function ${\mit \Omega}(u)$ the Xi function on the critical line ${\mit \Xi}(\I y)$ is rapidly decreasing in infinity. Therefore it is not easy to represent the real-valued function ${\mit \Xi}(\I y)$ with its rapid oscillations under the envelope of rapid decrease for increasing variable $y$ graphically in a large region of this variable $y$. An appropriate real amplification envelope is seen from (\ref{xizeta}) to be $\alpha(y)=\left|\frac{1}{\left(\frac{1+\I\,2y}{4}\right)!}\right|\frac{2\pi^{\frac{1}{4}}}{\sqrt{1+4y^2}}$ which rises ${\mit \Xi}(\I y)$ to the level of the Riemann zeta function $\zeta\left(\frac{1}{2}+\I t\right)$ on the critical line $z=\I y$. This is shown in Fig. 2.3.
\begin{figure}[h]
\includegraphics[width=16.0cm]{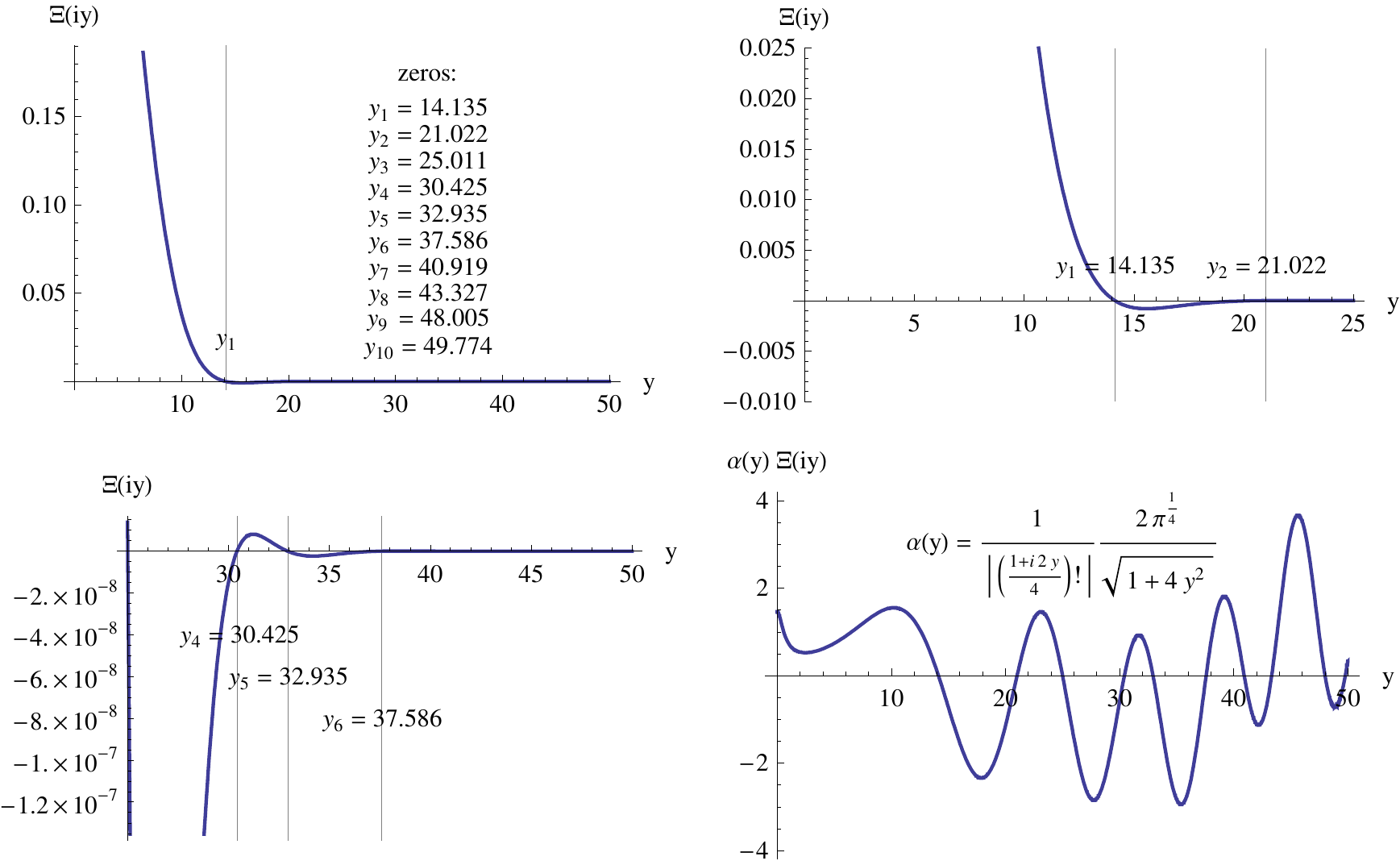}\hspace{1cm}
\caption{\foot Xi Function ${\mit \Xi}(\I y)$ on the imaginary axis $z=\I y$ (corresponding to $s=\frac{1}{2}+\I y$).
\newline \scriptsize{The envelope over the oscillations of the real-valued function ${\mit \Xi}(\I y)$ decreases extremely rapidly with increase of the variable $y$ in the shown intervals. This behavior makes it difficult to represent this function graphically for large intervals of the variable $y$. By an enhancement factor which rises the amplitude to the level of the zeta function $\zeta(s)$ we may see the oscillations under the envelope (last partial picture). A similar picture one obtains for the modulus of the Riemann zeta function $\left|\zeta\left(\frac{1}{2}+\I y\right)\right|$ only with our negative parts folded to the positive side of the abscissa, i.e. $\left|{\mit \Xi}(\I y)\right|=\left|\zeta\left(\frac{1}{2}+\I y\right)\right|$ (see also Fig. 2.1 (last partial picture)). The given values for the zeros at $\frac{1}{2}\pm \I y_n$ were first calculated by J.-P. Gram in 1903 up to $y_{15}$ \cite{edw}. We emphasize here that the shown very rapid decrease of the Xi function at the beginning of $y$ and for $y\rightarrow \pm \infty$ is due to the 'very high' smoothness of ${\mit \Omega}(u)$ for arbitrary $u$.}}
\end{figure}
The partial picture for $\alpha(y){\mit \Xi}(\I y)$ in Fig. 2.3. with negative part folded up is identical with the absolute value $\left|\zeta\left(\frac{1}{2}+\I t\right)\right|$ of the Riemann zeta function $\zeta(s)$ on the imaginary axis $s=\frac{1}{2}+\I t$ (fourth partial picture in Fig. 2.1).

We now give a representation of the Xi function by the derivative of the Omega function.
Using $\,{\rm ch}\left(uz\right)= \frac{1}{z}\frac{\partial}{\partial u}\,{\rm sh}\left(uz\right)$ one obtains from (\ref{capxirep3}) by partial integration the following alternative representation of the function ${\mit \Xi}\left(z\right)$
\begin{eqnarray}
{\mit \Xi}\left(z\right) &=& -\frac{1}{z}\int_{0}^{+\infty}du\,{\mit \Omega}^{(1)}(u)\,{\rm sh}\left(uz\right),\label{capxirep5}
\end{eqnarray}
that due to antisymmetry of $\,{\mit \Omega}^{(1)}(u)$ and $\,{\rm sh}\left(uz\right)$ with respect to $u\rightarrow -u$ can also be written
\begin{eqnarray}
{\mit \Xi}\left(z\right) &=& -\frac{1}{2z}\int_{-\infty}^{+\infty}du\,{\mit
\Omega}^{(1)}(u)\,{\rm sh}\left(uz\right) \;=\; -\frac{1}{2z}\int_{-\infty}^{+\infty}du\,{\mit
\Omega}^{(1)}\left(u\right)\,\E^{uz}. \label{capxirep6}
\end{eqnarray}
Figure 2.2 gives a graphical representation of the function ${\mit \Omega}\left(u\right)$ and of its first derivative ${\mit \Omega}^{(1)}(u) \equiv \frac{\partial {\mit \Omega}}{\partial u}(u)$ which due to rapid convergence of the sums is easily to generate by computer. One can express ${\mit \Xi}\left(z\right)$ also by higher derivatives ${\mit \Omega}^{(n)}(u) \equiv \frac{\partial^n {\mit \Omega}}{\partial u^n}(u)$ of the Omega function ${\mit \Omega}(u)$
according to
\begin{eqnarray}
{\mit \Xi}\left(z\right) &=& \frac{1}{z^{2m}}\int_{0}^{+\infty}du\,{\mit \Omega}^{(2m)}(u)\,{\rm ch}(uz) \nonumber\\ &=& -\frac{1}{z^{2m+1}}\int_{0}^{+\infty}du\,{\mit \Omega}^{(2m+1)}(u)\,{\rm sh}(uz),\quad \left(m=0,1,2,\ldots \right),\label{capxirep7}
\end{eqnarray}
with the symmetries of the derivatives of the function ${\mit \Omega}(u)$ for $u\leftrightarrow -u$
\begin{eqnarray}
{\mit \Omega}^{(2m)}(u) \;=\; +{\mit \Omega}^{(2m)}(-u),\quad {\mit \Omega}^{(2m+1)}(u) \;=\; -{\mit \Omega}^{(2m+1)}(-u),\;\; \Rightarrow \;\; {\mit \Omega}^{(2m+1)}(0) \;=\; 0,\quad \left(m=0,1,\ldots\right).\quad
\end{eqnarray}
This can be seen by successive partial integrations in (\ref{capxirep3}) together with complete induction. The functions ${\mit \Omega}^{(n)}(u)$ in these integral transformations are for $n \ge 1$ not monotonic functions.

We mention yet another representation of the function ${\mit \Xi}(z)$. Using the transformations
\begin{eqnarray}
t_n\equiv \pi n^2 \E^{2u}, \quad \Rightarrow \quad dt_n \;=\;2\pi n^2 \E^{2u}du \;=\; 2t_ndu,\quad du\;=\;\frac{dt_n}{2u},
\end{eqnarray}
the function ${\mit \Xi}(z)$ according to (\ref{capxirep4}) with the explicit representation of the function ${\mit \Omega}(u)$ in (\ref{defcapomega}) can now be represented in the form
\begin{eqnarray}
{\mit \Xi}(z) &=& \sum_{n=1}^{\infty}\frac{1}{(\pi n^2)^{\frac{1}{4}}}\left\{2\left(\left(\pi n^2\right)^{-\frac{\scr{z}}{2}}{\rm \Gamma}\left(\frac{9}{4}+\frac{z}{2},\pi n^2\right)+\left(\pi n^2\right)^{\frac{\scr{z}}{2}}{\rm \Gamma}\left(\frac{9}{4}-\frac{z}{2},\pi n^2\right)\right)\right. \nonumber\\&& \left. -3\left(\left(\pi n^2\right)^{-\frac{\scr{z}}{2}}{\rm \Gamma}\left(\frac{5}{4}+\frac{z}{2},\pi n^2\right)+\left(\pi n^2\right)^{\frac{\scr{z}}{2}}{\rm \Gamma}\left(\frac{5}{4}-\frac{z}{2},\pi n^2\right)\right)\right\},\label{xiincomplgamma}
\end{eqnarray}
where ${\rm \Gamma}(\alpha,x)$ denotes the incomplete Gamma function defined by (e.g., \cite{bate1,nist,paris}
\begin{eqnarray}
{\rm \Gamma}(\alpha,x) & \equiv & \int_{x}^{+\infty}dt\,\E^{-t}t^{\alpha-1} \;\equiv\; {\rm \Gamma}(\alpha)-{\rm \gamma}(\alpha,x).
\end{eqnarray}
However, we did not see a way to prove the Riemann hypothesis via the representation (\ref{xiincomplgamma}).

The Riemann hypothesis for the zeta function $\zeta(s = \sigma +\I t)$ is now equivalent to the hypothesis that all zeros of the related entire function ${\mit \Xi}\left(z=x+\I y\right)$ lie on the imaginary axis $z=\I y$ that means on the line to real part $x=0$ of $z=x+\I y$ which becomes now the critical line. Since the zeta function $\zeta(s)$ does not possess zeros in the convergence region $\sigma >1$ of the Euler product (\ref{eulerprod}) and due to symmetries (\ref{omegasym}) and (\ref{xisym}) it is only necessary to prove that ${\mit \Xi}(z)$ does not
possess zeros within the strips $-\frac{1}{2}\le x <0$ and $0 < x \le +\frac{1}{2}$ to both sides of the imaginary axis $z=\I y$ where for symmetry the proof for one of these strips would be already sufficient. However, we will go another way where the restriction to these strips does not play a role for the proof.

\setcounter{chapter}{3}
\setcounter{equation}{0}
\section*{3. Application of second mean-value theorem of calculus to Xi function}

After having accepted the basic integral representation (\ref{capxirep3}) of the entire function ${\mit \Xi}(z)$ according to
\begin{eqnarray}
{\mit \Xi}(z) & \equiv & \int_{0}^{+\infty}du\,{\mit \Omega}(u)\,{\rm ch}(uz),\label{basiccapxi0}
\end{eqnarray}
with the function ${\mit \Omega}(u)$ explicitly given in (\ref{defcapomega}) we concentrate us on its further treatment. However, we do this not with this specialization for the real-valued function ${\mit \Omega}(u)$ but with more general suppositions for it. Expressed by real part $U(x,y)$ and imaginary part $V(x,y)$ of ${\mit \Xi}(z)$
\begin{eqnarray}
{\mit \Xi}(x+\I y) \; \equiv \; U(x,y) +\I V(x,y), \quad U(x,y) \;=\; \left(U(x,y)\right)^*, \quad V(x,y) \;=\; \left(V(x,y)\right)^*,
\end{eqnarray}
we find from (\ref{basiccapxi0})
\begin{eqnarray}
U(x,y) \;=\; \int_{0}^{+\infty}du\,{\mit \Omega}(u)\,{\rm ch}(ux)\cos(uy), \quad V(x,y) \;=\; \int_{0}^{+\infty}du\,{\mit \Omega}(u)\,{\rm sh}(ux)\sin(uy).\label{reimpartxi}
\end{eqnarray}

We suppose now as necessary requirement for ${\mit \Omega}(u)$ and satisfied in the special case (\ref{defcapomega})
\begin{eqnarray}
{\mit \Omega}(u) \; \ge\; 0,\quad (0\le u < +\infty),\quad \Rightarrow \quad {\mit \Omega}(0) \; > \; 0.\label{condsecmvth}
\end{eqnarray}
Furthermore, ${\mit \Xi}(z)$ should be an entire function that requires that the integral (\ref{basiccapxi0}) is finite for arbitrary
complex $z$ and therefore that ${\mit \Omega}(u)$ is rapidly decreasing in infinity, more precisely
\begin{eqnarray}
\lim_{u\rightarrow +\infty}\frac{{\mit \Omega}(u)}{\exp\left(-\lambda u\right)} = 0, \quad 0 < 0 \le \lambda < +\infty ,\label{rapiddecr}
\end{eqnarray}
for arbitrary $\lambda \ge 0$. This means that the function ${\mit \Omega}(u)$ should be a nonsingular function which is rapidly decreasing in infinity, more rapidly than any exponential function $\E^{-\lambda u}$ with arbitrary $\lambda >0$. Clearly, this is satisfied for the special function ${\mit \Omega}(u)$ in (\ref{defcapomega}).

Our conjecture for a longer time was that all zeros of ${\mit \Xi}(z)$ lie on the imaginary axis $z=\I y$ for a large class of functions ${\mit \Omega}(u)$ and that this is not very specific for the special function ${\mit \Omega}(u)$ given in (\ref{defcapomega}) but is true for a much larger class. It seems that to this class belong all non-increasing functions ${\mit \Omega}(u)$, i.e such functions for which holds ${\mit \Omega}^{(1)}(u) \le 0$ for its first derivative and which rapidly decrease in infinity. This means that they vanish more rapidly in infinity than any power functions $|u|^{-n},\,(n=1,2,\ldots)$ (practically they vanish exponentially). However, for the convergence of the integral (\ref{basiccapxi0}) in the whole complex $z$-plane it is necessary that the functions have to decrease in infinity also more rapidly than any exponential function $\exp(-\lambda u)$ with arbitrary $\lambda >0$ expressed in (\ref{rapiddecr}). In particular, to this class belong all rapidly decreasing functions ${\mit \Omega}(u)$ which vanish from a certain $u\ge u_0$ on and which may be called non-increasing finite functions (or functions with compact support). On the other side, continuity of its derivatives ${\mit \Omega}^{(n)}(u),\,(n=1,2,\ldots)$ is not required. The modified Bessel functions ${\rm I}_{\nu}(z)$ 'normalized' to the form of entire functions $\left(\frac{2}{z}\right)^{\nu}{\rm I}_{\nu}(z)$ for $\nu \ge \frac{1}{2}$ possess a representation of the form (\ref{basiccapxi0}) with a function ${\mit \Omega}(u)$ which vanishes from $u \ge 1$ on but a number of derivatives of ${\mit \Omega}(u)$ for the functions is not continuous at $u=1$ depending on the index $\nu$. It is valuable that here an independent proof of the property that all zeros of the modified Bessel functions ${\rm I}_{\nu}(u)$ lie on the imaginary axis can be made using their differential equations via duality relations. We intend to present this in detail in a later work.

Furthermore, to the considered class belong all monotonically  decreasing functions with the described rapid decrease in infinity. The fine difference of the decreasing functions to the non-increasing functions ${\mit \Omega}(u)$ is that in first case the function ${\mit \Omega}(u)$ cannot stay on the same level in a certain interval that means we have ${\mit \Omega}^{(1)}(u) <0$ for all points $u > 0$ instead of ${\mit \Omega}^{(1)}(u) \le 0$ only. A function which decreases not faster than $\E^{-\lambda u}$ in infinity does not fall into this category as, for example, the function ${\rm sech}(z) \equiv \frac{1}{{\rm ch}(z)}$ shows. On the other side, also some simply calculable discrete superpositions such as $a_1\,{\rm ch}(u)+a_2\,{\rm ch}(2u)$ or $a_1\,{\rm ch}(z) +a_3\,{\rm ch}(3z)$ as function ${\mit \Xi}(z)$ with positive amplitudes $a_n$ do not provide a counterexample that the zeros lie outside the imaginary axis but show that if the amplitudes $a_n$ do not possess a definite sign then they may possess zeros outside the imaginary axis.

To apply the second mean-value theorem it is necessary to restrict us to a class of functions ${\mit \Omega}(u) \rightarrow f(u)$ which are non-increasing that means for which for all $u_1 < u_2$ in considered interval holds
\begin{eqnarray}
f(a) \ge f(u_1) \ge f(u_2) \ge f(b) \ge 0,\quad (a \le u_1 \le u_2 \le b),
\end{eqnarray}
or equivalently in more compact form
\begin{eqnarray}
f^{(1)}(u) & \le & 0, \quad (a \le u \le b).
\end{eqnarray}
In case of $f^{(1)}(u)=0$ for certain $u$ the next higher non-vanishing derivative should be negative.
The monotonically  decreasing functions in the interval $a\le u \le b$, in particular, belong to the class of non-increasing functions with the fine difference that here
\begin{eqnarray}
f(a) > f(u_1) > f(u_2) > f(b) > 0,\quad (a < u_1 < u_2 < b),
\end{eqnarray}
is satisfied. If furthermore $g(u)$ is a continuous function in the interval $a \le u \le b$ the second mean-value theorem (often called theorem of Bonnet (1867) or Gauss-Bonnet theorem) states an equivalence for the following integral on the left-hand side to the expression on the right-hand side according to (see some monographs about calculus or real analysis; we recommend the monographs of Courant \cite{courant} (Appendix to chap IV) and of Widder \cite{widder} who called it Weierstrass form of Bonnet's theorem (chap. 5, \S 4))
\begin{eqnarray}
\int_{a}^{b}du\,f(u)g(u) &=& f(a)\int_{a}^{u_0}du\,g(u) +f(b)\int_{u_0}^{b}du\,g(u),\quad (a\le u_0 \le b),\label{secmeanvth}
\end{eqnarray}
where $u_0 $ is a certain value within the interval boundaries $a < b $ which as a rule we do not exactly know. It holds also for non-decreasing functions which include the monotonically  increasing functions as special class in analogous way. The proof of the second mean-value theorem is comparatively simple by applying a substitution in the (first) mean-value theorem of integral calculus \cite{courant,widder}.

Applied to our function $f(u) = {\mit \Omega}(u)$ which in addition should rapidly decrease in infinity according to (\ref{rapiddecr})  this means in connection with monotonic decrease that it has to be positively semi-definite if ${\mit \Omega}(0) > 0$ and therefore
\begin{eqnarray}
{\mit \Omega}(0) \ge {\mit \Omega}(u) \ge 0,\quad {\mit \Omega}^{(1)}(u) \le 0,\quad \left(0\le u \le +\infty\right),\quad
{\mit \Omega}(u\rightarrow +\infty) \rightarrow 0,\label{secmeanvth0}
\end{eqnarray}
and the theorem (\ref{secmeanvth}) takes on the form
\begin{eqnarray}
\int_{0}^{+\infty}du\,{\mit \Omega}(u) g(u) &=& {\mit \Omega}(0)\int_{0}^{u_0}du\,g(u),\quad (0 \le u_0 < +\infty), \label{secmeanvth1}
\end{eqnarray}
where the extension to an upper boundary $b\rightarrow +\infty$ in (\ref{secmeanvth}) for $f(+\infty)=0$ and in case of existence of the integral is unproblematic.

If we insert in (\ref{secmeanvth}) for $g(u)$ the function ${\rm ch}(uz)$ which apart from the real variable $u$ depends in parametrical way on the complex variable $z$ and is an analytic function of $z$ we find that $u_0$ depends on this complex parameter also in an analytic way as follows
\begin{eqnarray}
{\mit \Xi}(z) & \equiv & \int_{0}^{+\infty}du\,{\mit \Omega}(u)\,{\rm ch}(uz) \;=\; {\mit \Omega}(0)\int_{0}^{w_0(z)}du\,{\rm ch}(uz) \nonumber\\ &=& {\mit \Omega}(0)\frac{\,{\rm sh}\left(w_0(z)z\right)}{z}, \quad w_0(z) \;=\; u_0(x,y) +\I v_0(x,y),\label{secmeanvth2compl}
\end{eqnarray}
where $w_0(x+\I z)=u_0(x,y)+\I v_0(x,y)$ is an entire function with $u_0(x,y)$ its real and $v_0(x,y)$ its imaginary part. The condition for zeros $z\neq 0$ is that ${\rm sh}\left(w_0(z)z\right){z}$ vanishes that leads to
\begin{eqnarray}
w_0(z)z &=& \left(u_0(x,y)+\I v_0(x,y)\right)\left(x+\I y\right) \;=\; \I n\pi, \quad (n=0,\pm 1,\pm 2,\ldots), \label{zcuvxy0}
\end{eqnarray}
or split in real and imaginary part
\begin{eqnarray}
u_0(x,y)x -v_0(x,y)y \;=\; 0,\label{zcuvxyre}
\end{eqnarray}
for the real part and
\begin{eqnarray}
u_0(x,y)y +v_0(x,y)x \;=\; n\pi, \quad (n=0,\pm 1,\pm 2,\ldots),\label{zcuvxyim}
\end{eqnarray}
for the imaginary part.

The multi-valuedness of the mean-value functions in the conditions (\ref{zcuvxy0}) or (\ref{zcuvxyim}) is an interesting phenomenon which is connected with the periodicity of the function $g(u)=\,{\rm ch}(uz)$ on the imaginary axis $z=\I y$ in our application (\ref{secmeanvth2compl}) of the second mean-value theorem (\ref{secmeanvth1}). To our knowledge this is up to now not well studied. We come back to this in the next Sections 4 and, in particular, Section 7 brings some illustrative clarity when we represent the mean-value functions graphically. At present we will say only that we can choose an arbitrary $n$ in (\ref{zcuvxyim}) which provides us the whole spectrum of zeros $z_1,z_2,\ldots$  on the upper half-plane and the corresponding spectrum of zeros $z_{-1}=-z_1,z_{-2}=-z_2,\ldots $ on the lower half-plane of ${\mathbb C}$ which as will be later seen lie all on the imaginary axis. Since in computer calculations the values of the Arcus Sine function are provided in the region from $-\frac{\pi}{2}$ to $+\frac{\pi}{2}$ it is convenient to choose $n=0$ but all other values of $n$ in (\ref{zcuvxyim}) lead to equivalent results.

One may represent the conditions (\ref{zcuvxyre}) and (\ref{zcuvxyim}) also in the following equivalent form
\begin{eqnarray}
u_0(x,y)\;=\; \frac{y}{x^2+y^2}n\pi,\quad v_0(x,y)\;=\; \frac{x}{x^2+y^2}n\pi,\label{zcuvxy2}
\end{eqnarray}
from which follows
\begin{eqnarray}
\left(u^2_0(x,y)+v^2_0(x,y)\right)\left(x^2+y^2\right) \;=\;(n\pi)^2,\quad \frac{v_0(x,y)}{u_0(x,y)} \;=\; \frac{x}{y}.\label{zcuvxy3}
\end{eqnarray}

All these forms (\ref{zcuvxyre})--(\ref{zcuvxy3}) are implicit equations with two variables $(x,y)$ which cannot be resolved with respect to one variable (e.g., in forms $y=y_k(x)$ for each fixed $n$  and branches $k$) and do not provide immediately the necessary conditions for zeros in explicit form but we can check that (\ref{zcuvxy2}) satisfies the Cauchy-Riemann equations as a minimum requirement
\begin{eqnarray}
\frac{\partial u_0(x,y)}{\partial x} \;=\; \frac{\partial v_0(x,y)}{\partial y},\quad \frac{\partial u_0(x,y)}{\partial y} \;=\; -\frac{\partial v_0(x,y)}{\partial x}.
\end{eqnarray}
We have to establish now closer relations between real and imaginary part $u_0(x,y)$ and $v_0(x,y)$ of the complex mean-value parameter $w_0(z=x+\I y)$. The first step in preparation to this aim is the consideration of the derived conditions on the imaginary axis.

\setcounter{chapter}{4}
\setcounter{equation}{0}
\section*{4. Specialization of second mean-value theorem to Xi function on imaginary axis}

By restriction to the real axis $y=0$ we find from (\ref{reimpartxi}) for the function ${\mit \Xi}(z)$
\begin{eqnarray}
{\mit \Xi}(x) &=& U(x,0), \quad  V(x,0) \;=\; 0,
\end{eqnarray}
with the following two possible representations of $U(x,0)$ related by partial integration
\begin{eqnarray}
U(x,0) &=& \int_{0}^{+\infty}du\,{\mit \Omega}(u)\,{\rm ch}(ux) \;=\; -\frac{1}{x} \int_{0}^{+\infty}du\,{\mit \Omega}^{(1)}(u)\,{\rm sh}(ux) \; > \; 0.\label{reimpartofimxi}
\end{eqnarray}
The inequality $U(x,0)>0$ follows according to the supposition ${\mit \Omega}(u) \ge 0, {\mit \Omega}(0) >0$ from the non-negativity of the integrand that means from ${\mit \Omega}(u)\,{\rm ch}(ux) \ge 0$. Therefore, the case $y=0$ can be excluded from the beginning in the further considerations for zeros of $U(x,y)$ and $V(x,y)$.

We now restrict us to the imaginary axis $x=0$ and find from (\ref{reimpartxi}) for the function ${\mit \Xi}(z)$
\begin{eqnarray}
{\mit \Xi}(\I y) &=& U(0,y), \quad  V(0,y) \;=\; 0.
\end{eqnarray}
with the following two possible representations of $U(0,y)$ related by partial integration
\begin{eqnarray}
U(0,y) &=& \int_{0}^{+\infty}du\,{\mit \Omega}(u)\cos(uy) \;=\; -\frac{1}{y} \int_{0}^{+\infty}du\,{\mit \Omega}^{(1)}(u)\sin(uy).\label{reimpartofimxi1}
\end{eqnarray}

From the obvious inequality
\begin{eqnarray}
-1 \; \le \; \cos(uy) \; \le +1,
\end{eqnarray}
together with the supposed positivity of ${\mit \Omega}(u)$ one derives from the first representation of $U(0,y)$ in (\ref{reimpartofimxi1}) the inequality
\begin{eqnarray}
-{\mit \Omega}_0 \; \le\; U(0,y) \; \le\;+{\mit \Omega}_0,\quad {\mit \Omega}_0 \; =\; U(0,0) \;\equiv \; \int_{0}^{+\infty}du\,{\mit \Omega}(u) \; \ge \; 0.\label{rexiineq1}
\end{eqnarray}
In the same way by the inequality
\begin{eqnarray}
-1 \; \le \; \sin(uy) \; \le +1,
\end{eqnarray}
one derives using the non-positivity of ${\mit \Omega}^{(1)}(u)$ (see (\ref{secmeanvth0})) together with the second representation of $U(0,y)$ in $(\ref{reimpartofimxi1})$ the inequality
\begin{eqnarray}
-{\mit \Omega}(0) \; \le \; U(0,y)y \; \le \; +{\mit \Omega}(0),\quad {\mit \Omega}(0) \;=\; -\int_{0}^{+\infty}du\,{\mit \Omega}^{(1)}(u) \; \ge \; 0.\label{rexiineq2}
\end{eqnarray}
which as it is easily seen does not depend on the sign of $y$. Therefore we have two non-negative parameters, the zeroth moment ${\mit \Omega}_0$ and the value ${\mit \Omega}(0)$, which according to (\ref{rexiineq1}) and (\ref{rexiineq2}) restrict the range of values of $U(0,y)$ to an interior range both to (\ref{rexiineq1}) and to (\ref{rexiineq2}) at once.

For mentioned purpose we now consider the restriction of the mean-value parameter $w_0(z)$ to the imaginary axis $z=\I y$ for which $g(u) = {\rm ch}(u(\I y))=\cos(uy)$ is a real-valued function of $y$. For arbitrary fixed $y$ we find by the second mean-value theorem a parameter $u_0$ in the interval $0 \le y < +\infty$ which naturally depends on the chosen value $y$ that means $u_0=u_0(0,y)$. The extension from the imaginary axis $z=\I y$ to the whole complex plane ${\mathbb C}$ can be made then using methods of complex analysis. We discuss some formal approaches to this in Appendix B. Now we apply (\ref{secmeanvth2compl}) to the imaginary axis $z=\I y$.

The second mean-value theorem (\ref{secmeanvth2compl}) on the imaginary axis $z=\I y$ (or $x=0$) takes on the form
\begin{eqnarray}
\int_{0}^{+\infty}du\,{\mit \Omega}(u)\,{\rm ch}(u(\I y)) &=& \int_{0}^{+\infty}du\,{\mit \Omega}(u)\cos(u y) \;=\; {\mit \Omega}(0) \int_{0}^{u_0(0,y)}du\,\cos(uy) \nonumber\\ &=& {\mit \Omega}(0)\frac{\sin\left(u_0(0,y) y\right)}{y},\quad \left(u_0(0,y) \neq  0,\quad v_0(0,y) =0 \right).
\label{realsecmeanvth}
\end{eqnarray}
As already said since the left-hand side is a real-valued function the right-hand side has also to be real-valued and the parameter function $w_0(\I y)$ is real-valued and therefore it can only be the real part $u_0(0,y)$ of the complex function $w_0(z= x +\I y)= u_0(x,y) +\I v_0(x,y)$ for $x=0$.

The second mean-value theorem states that $u_0(0,y)$ lies between the minimal and maximal values of the integration borders that is here between $0$ and $+\infty$ and this means that $u_0(0,y)$ should be positive. Here arises a problem which is connected with the periodicity of the function $g(u)=\cos(uy)$ as function of the variable $u$ for fixed variable $y$ in the application of the mean-value theorem. Let us first consider the special case $y=0$ in (\ref{realsecmeanvth}) which leads to
\begin{eqnarray}
\int_{0}^{+\infty}du\,{\mit \Omega}(u) &=& {\mit \Omega}(0) \int_{0}^{u_0}du \;=\; {\mit \Omega}(0)u_0 \;=\;
{\mit \Omega}(0)u_0\underbrace{\lim_{y\to 0}\frac{\sin(u_0y)}{u_0y}}_{=\,1},\quad u_0\,\equiv \; u_0(0,0)\;>\;0.
\end{eqnarray}
From this relation follows $u_0 \equiv u_0(0,0)>0$ and it seems that all is correct also with the continuation to $u_0(0,y) > 0$ for arbitrary $y$. One may even give the approximate values ${\mit \Omega}(0) \approx 1.78679 $ and $u_0 \approx 0.27822 $ and therefore ${\mit \Omega}_0 \equiv {\mit \Omega}(0)u_0 \approx 0.49712 $ which, however, are not of importance for the later proofs. If we now start from $u_0(0,0) > 0$ and continue it continuously to $u_0(0,y)$ then we see that $u_0(0,y)$ goes monotonically to zero and approaches zero approximately at $y=y_1\approx 14.135$ that is at the first zero of the function ${\mit \Xi}(\I y)$ on the positive imaginary axis and goes then first beyond zero and oscillates then with decreasing amplitude for increasing $y$ around the value zero with intersecting it exactly at the zeros of ${\mit \Xi}(\I y)$. We try to illustrate this graphically in Section 7. All zeros lie then on the branch $u_0(0,y)y=n\pi$ with $n=0$. That $u_0(0,y)$ goes beyond zero seems to contradict the content of the second mean-value theorem according which $u_0(0,y)$ has to be positive in our application. Here comes into play the multi-valuedness of the mean-value function $u_0(0,y)$. For the zeros of $\sin(u_0(0,y)y)$ in (\ref{realsecmeanvth}) the relations $u_0(0,y)y =n\pi$ with different integers $n$ are equivalent  and one may find to values $u_0(0,y) <0$ equivalent curves $u_0(n;0,y)$ with $u_0(n;0,y)>0$ and all these curves begin with $u_0(n\neq 0;0,0) \rightarrow \infty$ for $y\rightarrow 0$. However, we cannot continue $u_0(0,0)$ in continuous way to only positive values for $u_0(0,y)$.

For $|y| \rightarrow \infty$ the inequality (\ref{rexiineq2}) is stronger than (\ref{rexiineq1}) and characterizes the restrictions of $U(0,y)$ and via the equivalence $U(0,y)y= {\mit \Omega}(0)\sin(u_0(0,y)y)$ follows from (\ref{rexiineq2})
\begin{eqnarray}
\left(n-\frac{1}{2}\right)\pi \;\le \; u_0(0,y)y \;=\; \arcsin\left(\frac{U(0,y)y}{{\mit \Omega}(0)}\right) \; \le \; \left(n+\frac{1}{2}\right)\pi,\quad (n=0,\pm 1, \pm 2,\ldots),\label{basicinterval}
\end{eqnarray}
where the choice of $n$ determines a basis interval of the involved multi-valued function $\arcsin(z)$ and the inequality says that it is in every case possible to choose it from the same interval of length $\pi$. The zeros $y_k$ of the Xi function ${\mit \Xi}(x+\I y)$ on the imaginary axis $x=0$ (critical line) are determined alone by the (multi-valued) function $u_0(0,y)$ whereas $v_0(0,y)$ vanishes automatically on the imaginary axis in considered special case and does not add a second condition. Therefore, the zeros are the solutions of the conditions
\begin{eqnarray}
u_0(0,y)y \;=\; n\pi,\quad (n=0,\pm 1, \pm2, \ldots),\quad (v_0(0,y) \;=\;0).\label{zerosimax}
\end{eqnarray}
It is, in general, not possible to obtain the zeros $y_k$ on the critical line exactly from the mean-value function $u_0(0,y)$ in (\ref{realsecmeanvth}) since generally we do not possess it explicitly.

In special cases the function $u_0(0,y)$ can be calculated explicitly that is the case, for example, for all (modified) Bessel functions $\left(\frac{2}{z}\right)^{\nu}{\rm I}_{\nu}(z)$. The most simple case among these is the case $\nu=\frac{1}{2}$ when the corresponding function ${\mit \Omega}(u)$ is a step function
\begin{eqnarray}
&& {\mit \Omega}(u) \;=\; {\mit \Omega}(0)\theta(u_0-u),\label{omegastepf}
\end{eqnarray}
where $\theta(x) = \left\{\begin{array}{cc} 0,& x<0 \\ 1,& x>0 \end{array}\right.$ is the Heaviside step function. In this case follows
\begin{eqnarray}
&& {\mit \Xi}(z) \;=\; {\mit \Omega}(0)\int_{0}^{u_0}du\,{\rm ch}(uz) \;=\; {\mit \Omega}(0)u_0\, \frac{{\rm sh}(u_0 z)}{u_0z} \;=\; \underbrace{{\mit \Omega}(0)u_0}_{=\,{\mit \Omega}_0}\, \left(\frac{\pi}{2u_0z}\right)^{\frac{1}{2}}\,{\rm I}_{\frac{1}{2}}(u_0z),\label{stepfunckernel}
\end{eqnarray}
where ${\mit \Omega}(0)u_0 =\int_{0}^{+\infty}du\,{\mit \Omega}(u)$ is the area under the function ${\mit \Omega}(u)= {\mit \Omega}(0)\theta(u_0-u)$ (or the zeroth-order moment of this function. For the squared modulus of the function ${\mit \Xi}(z) $ we find
\begin{eqnarray}
{\mit \Xi}(z)\left({\mit \Xi}(z)\right)^* \;=\; \left({\mit \Omega}(0)\right)^2\frac{{\rm sh}^2(u_0x)+\sin^2(u_0y)}{x^2+y^2}
\;=\; \left({\mit \Omega}(0)\right)^2 \frac{{\rm ch}(2u_0x)-\cos(2u_0y)}{2\left(x^2+y^2\right)},\label{xixistarstep}
\end{eqnarray}
from which, in particular, it is easy to see that this special function ${\mit \Xi}(x+\I y)$ possesses zeros only on the imaginary axis $z=\I y$ or $x=0$ and that they are determined by
\begin{eqnarray}
u_0y_n &=& n\pi, \quad \Rightarrow \quad y_n \;=\; \frac{n\pi}{u_0},\quad (n=\pm 1,\pm 2,\ldots).
\end{eqnarray}
The zeros on the imaginary axis are here equidistant but the solution $y_0=0$ is absent since then also the denominators in (\ref{xixistarstep}) are vanishing. The parameter $w_0(z)$ in the second mean-value theorem is here a real constant $u_0$ in the whole complex plane
\begin{eqnarray}
&& w_0(x+\I y) \;=\; u_0(x,y)+\I v_0(x,y) \;=\;u_0,\quad \Rightarrow \nonumber\\ && u_0(x,y) \;=\; u(0,y) \;=\; u_0,\quad v_0(x,y)\;=\;v_0(0,y) \;=\;0.\label{stepfunckernel1}
\end{eqnarray}
Practically, the second mean-value theorem compares the result for an arbitrary function ${\mit \Omega}(u)$ under the given restrictions with that for a step function ${\mit \Omega}(u)={\mit \Omega}(0)\,\theta(u_0-u)$ by preserving the value ${\mit \Omega}(0)$ and making the parameter $u_0$ depending on $z$ in the whole complex plane. Without discussing now quantitative relations the formulae (\ref{stepfunckernel1}) suggest that $v_0(x,y)$ will stay a 'small' function compared with $u_0(x,y)$ in the neighborhood of the imaginary axis (i.e. for $|x|\ll |y|$) in a certain sense.

We will see in next Section that the function $u_0(0,y)$ taking into account $v_0(0,y)=0$ determines the functions $u_0(x,y)$ and $v_0(x,y)$ and thus $w_0(z)$ in the whole complex plane via the Cauchy-Riemann equations in an operational approach that means in an integrated form which we did not found up to now in literature. The general formal part is again delegated to an Appendix B.

\setcounter{chapter}{5}
\setcounter{equation}{0}
\section*{5. Accomplishment of proof for zeros of Xi functions on imaginary axis alone}

In last Section we discussed the application of the second mean-value theorem to the function ${\mit \Xi}(z)$ on the imaginary axis $z=\I y$. Equations (\ref{zcuvxyre}) and (\ref{zcuvxyim})
or their equivalent forms (\ref{zcuvxy2}) or (\ref{zcuvxy3}) are not yet sufficient to derive conclusions about the position of the zeros on the imaginary axis in dependence on $x\neq 0$. We have yet to derive more information about the mean-value functions $w_0(z)$ which we obtain by relating the real-valued function $u_0(x,y)$ and $v_0(x,y)$ to the function $u_0(0,y)$ on the imaginary axis ($v_0(0,y)=0$).

The general case of complex $z$ can be obtained from the special case $z=\I y$ in (\ref{realsecmeanvth}) by application of the displacement operator $\exp\left(-\I x \frac{\partial}{\partial y}\right)$ to the function ${\mit \Xi}(\I y)$ according to
\begin{eqnarray}
{\mit \Xi}(x+\I y) &=& \exp\left(-\I x\frac{\partial}{\partial y}\right){\mit \Xi}(\I y) \exp\left(\I x\frac{\partial}{\partial y}\right) \;=\; \exp\left(-\I x\frac{\partial}{\partial y}\right)\int_{0}^{+\infty}du\,{\mit \Omega}(u)\,{\rm ch}(\I uy) \exp\left(\I x\frac{\partial}{\partial y}\right) \nonumber\\ &=& \exp\left(-\I x\frac{\partial}{\partial y}\right){\mit \Omega}(0)\frac{\,{\rm sh}\Big(\I u_0(0,y)y\Big)}{\I y}\exp\left(\I x\frac{\partial}{\partial y}\right)\;=\;
{\mit \Omega}(0)\frac{\,{\rm sh}\Big(u_0(0,y-\I x)\I (y-\I x)\Big)}{\I (y-\I x)}.\qquad \label{xixy}
\end{eqnarray}
The function $u_0(0,y-\I x) = w_0(x+\I y) = u_0(x,y) +\I v_0(x,y)$ is related to $u_0(0,y)$ as follows
\begin{eqnarray}
&& u_0(x,y) \;=\;
\cos\left(x\frac{\partial}{\partial y}\right)u_0(0,y), \quad v_0(x,y) \;=\;
-\sin\left(x\frac{\partial}{\partial y}\right)u_0(0,y), \label{uvxycossin}
\end{eqnarray}
or in more compact form
\begin{eqnarray}
w_0(x+\I y) &=& \exp\left(-\I x\frac{\partial}{\partial y}\right)u_0(0,y) \;=\; u_0(0,y-\I x).\label{uvxycossin1}
\end{eqnarray}
This is presented in Appendix B in more general form for additionally non-vanishing $v_0(0,y)$ and arbitrary holomorphic functions. It means that we may obtain $u_0(x,y) $ and $v_0(x,y)$ by applying the operators $\cos\left(x\frac{\partial}{\partial y}\right)$ and $-\sin\left(x\frac{\partial}{\partial y}\right)$, respectively, to the function $u_0(0,y)$ on the imaginary axis (remind $v_0(0,y)=0$ vanishes there in our case). Clearly, equations (\ref{uvxycossin}) are in agreement with the Cauchy-Riemann equations $\frac{\partial u_0}{\partial x}= \frac{\partial v_0}{\partial y}$ and $\frac{\partial u_0}{\partial y}= -\frac{\partial v_0}{\partial x}$ as a minimal requirement.

We now write ${\mit \Xi}(z)$ in the form equivalent to (\ref{xixy})
\begin{eqnarray}
{\mit \Xi}(x+ \I y) &=& {\mit \Omega}(0) \frac{{\rm sh}\Big(\big(u_0(x,y)+\I v_0(x,y)\big)(x+\I y)\Big)}{x+ \I y} \nonumber\\ &=& {\mit \Omega}(0) \frac{{\rm sh}\Big(\big(u_0(x,y)x- v_0(x,y)y\big)+\I\big(u_0(x,y)y + v_0(x,y)x\big)\Big)}{x+\I y}.\label{secmeanvth5}
\end{eqnarray}
The denominator $x+\I y$ does not contribute to zeros. Since the Hyperbolic Sine possesses zeros only on the imaginary axis we see from (\ref{secmeanvth5}) that we may expect zeros only for such related variables $(x,y)$ which satisfy the necessary condition of vanishing of its real part of the argument that leads as we already know to (see (\ref{zcuvxyre}))
\begin{eqnarray}
u_0(x,y)x- v_0(x,y)y &=& 0. \label{neccondzeros}
\end{eqnarray}
The zeros with coordinates $(x_k,y_k)$ themselves can be found then as the (in general non-degenerate) solutions of the following equation (see (\ref{zcuvxyim}))
\begin{eqnarray}
u_0(x,y)y + v_0(x,y)x &=& n\pi,\quad (n=0,\pm 1,\pm 2,\ldots), \label{zerosequ}
\end{eqnarray}
if these pairs $(x,y)$ satisfy the necessary condition (\ref{neccondzeros}). Later we will see provides the whole spectrum of solutions for the zeros but we can also obtain each $(x_k,y_k)$ separately from one branch $n$ and would they then denote by $(x_n,y_n)$.
Thus we have first of all to look for such pairs $(x,y)$ which satisfy the condition (\ref{neccondzeros}) off the imaginary axis that is for $x\neq 0$ since we know already that these functions may possess zeros on the imaginary axis $z=\I y$.

Using (\ref{uvxycossin}) we may represent the necessary condition (\ref{neccondzeros}) for the proof by the second mean-value theorem in the form
\begin{eqnarray}
x\cos\left(x\frac{\partial}{\partial y}\right)u_0(0,y) +y\sin\left(x\frac{\partial}{\partial y}\right)u_0(0,y) &=& 0,\label{zeroscondnec}
\end{eqnarray}
and equation (\ref{zerosequ}) which determines then the position of the zeros can be written with equivalent values $n$
\begin{eqnarray}
y\cos\left(x\frac{\partial}{\partial y}\right)u_0(0,y) - x\sin\left(x\frac{\partial}{\partial y}\right)u_0(0,y) &=& n\pi,\quad (n=0,\pm 1,\pm 2,\ldots). \label{zerosequ0}
\end{eqnarray}
We may represent Eqs. (\ref{zeroscondnec}) and (\ref{zerosequ0}) in a simpler form using the following operational identities
\begin{eqnarray}
x\cos\left(x\frac{\partial}{\partial y}\right)+y\sin\left(x\frac{\partial}{\partial y}\right) &=& \sin\left(x\frac{\partial}{\partial y}\right)y,\nonumber\\ y\cos\left(x\frac{\partial}{\partial y}\right)-x\sin\left(x\frac{\partial}{\partial y}\right) &=&  \cos\left(x\frac{\partial}{\partial y}\right)y,
\end{eqnarray}
which are a specialization of the operational identities (\ref{fydisplfxyop}) in Appendix B with $w(z) =u(x,y)+\I v(x,y) \rightarrow z=x+\I y $ and therefore $u(x,y)\rightarrow x,\,v(x,y) \rightarrow y$.
If we multiply (\ref{zeroscondnec}) and (\ref{zerosequ0}) both by the function $u_0(0,y)$ then we may write (\ref{zeroscondnec}) in the form (changing order $yu_0(0,y)=u_0(0,y)y$)
\begin{eqnarray}
\sin\left(x\frac{\partial}{\partial y}\right)\left(u_0(0,y)y\right) &=& 0,\label{zeroscond1}
\end{eqnarray}
and (\ref{zerosequ0}) in the form
\begin{eqnarray}
\cos\left(x\frac{\partial}{\partial y}\right)\left(u_0(0,y)y\right) &=& n\pi,\quad (n=0,\pm 1,\pm 2,\ldots). \label{zeroscond2}
\end{eqnarray}
The left-hand side of these conditions possess the general form for the extension of a holomorphic function $W(z)=U(x,y)+\I V(x,y)$ from the functions $U(0,y)$ and $V(0,y)$ on the imaginary axis to the whole complex plane in case of $V(0,y)=0$ and if we apply this to the function $U(0,y)=u_0(0,y)y$. Equations (\ref{zeroscond1}) and (\ref{zeroscond2}) possess now the most simple form, we found, to accomplish the proof for the exclusive position of zeros on the imaginary axis. All information about the zeros of the Xi function ${\mit \Xi}(z) = U(x,y)+\I V(x,y)$ for arbitrary $x$ is now contained in the conditions (\ref{zeroscond1}) and (\ref{zeroscond2}}) which we now discuss.

Since $\cos\left(x\frac{\partial}{\partial y}\right)$ is a nonsingular operator we can multiply both sides of equation (\ref{zeroscond2}) by the inverse operator $\cos^{-1}\left(x\frac{\partial}{\partial y}\right)$ and obtain
\begin{eqnarray}
u_0(0,y)y &=& \cos^{-1}\left(x\frac{\partial}{\partial y}\right) n\pi \;=\; \left\{1+\frac{1}{2}\left(x\frac{\partial}{\partial y}\right)^2 +\ldots\right\}n\pi \;=\;n\pi,\quad (n=0,\pm 1,\ldots).\qquad\label{zeroimaxcond1}
\end{eqnarray}
This equation is yet fully equivalent to (\ref{zeroscond2}) for arbitrary $x$ but it provides only the same solutions for the values $y$ of zeros as for zeros on the imaginary axis. This alone already suggests that it cannot be that zeros with $x\neq 0$ if they exist possess the same values of $y$ as the zeros on the imaginary axis. But in such form the proof of the impossibility of zeros off the imaginary axis seemed to be not satisfactory and we present in the following some slightly different variants which go deeper into the details of the proof.

In analogous way by multiplication of (\ref{zeroscond1}) with the operator $\sin\left(x\frac{\partial}{\partial y}\right)$ and (\ref{zeroscond2}) with the operator $\cos\left(x\frac{\partial}{\partial y}\right)$ and addition of both equations we also obtain condition (\ref{zeroimaxcond1}) that means
\begin{eqnarray}
u_0(0,y)y &=& \left\{\sin^2\left(x\frac{\partial}{\partial y}\right)+\cos^2\left(x\frac{\partial}{\partial y}\right)\right\}u_0(0,y)y \;=\; n\pi,\quad (n=0,\pm 1,\pm 2,\ldots),\label{zeroimaxcond2}
\end{eqnarray}
The equal conditions (\ref{zeroimaxcond1}) and (\ref{zeroimaxcond2}) which are identical with the condition for zeros on the imaginary axis are a necessary condition for all zeros. For each chosen equivalent $n$ (remind $u_0(0,y)$ depends then on $n$ which we do not mention by the notation) one obtains an infinite series of solutions $y_k$ for the zeros of the function ${\mit \Xi}(\I y)$
\begin{eqnarray}
u_0(0,y_k)y_k &=& n\pi, \qquad \{u_0(0,y)y\}_{y\neq y_k} \;\neq\; n\pi, \label{yk0}
\end{eqnarray}
whereas for $y\neq y_k$ equation (\ref{zeroimaxcond1}), by definition of $y_k$, is not satisfied. Supposing that we know $u_0(0,y)$ that is as a rule not the case, we could solve for each $n = 0,\pm 1, \pm 2,\ldots$ the usually transcendental equations (\ref{zeroimaxcond2}) graphically, for example, by drawing the equivalent functions $u_0(0,y)y$ over
variable $y$ as abscissa and looking for the intersections points with the lines $n\pi$ over $y$ (Section 7). These intersection points $y=y_k$ are the solutions for zeros $y_k$ on the imaginary axis. Choosing $x=0$ the second condition (\ref{zeroscond1}) is identically satisfied.

Now we have to look for additional zeros $(x,y_k)$ with $x\neq 0$. Whereas for zeros with $x=0$ the condition (\ref{zeroscond1}) is identically satisfied we have to examine this condition for zeros with $x\neq 0$. In the case of $x\neq 0$ we may divide both sides of the condition (\ref{zeroscond1}) by $x$ and obtain
\begin{eqnarray}
\frac{\sin\left(x\frac{\partial}{\partial y}\right)}{x}\left(u_0(0,y)y\right) &=& \frac{\sin\left(x\frac{\partial}{\partial y}\right)}{x\frac{\partial}{\partial y}} \frac{\partial}{\partial y}\left(u_0(0,y)y\right) \;=\; 0.\label{zeroscond3}
\end{eqnarray}
Since $\frac{\sin\left(x\frac{\partial}{\partial y}\right)}{x\frac{\partial}{\partial y}} $ is a nonsingular operator (in contrast to $\sin\left(x\frac{\partial}{\partial y}\right)$ which possesses $0$ as eigenvalue to eigenfunction $f(x,y)={\rm const}x^n,(n=0,1,\ldots)$) we may multiply equation (\ref{zeroscond3})
by the inverse operator $\frac{x\frac{\partial}{\partial y}}{\sin\left(x\frac{\partial}{\partial y}\right)}$ and obtain
\begin{eqnarray}
\frac{\partial}{\partial y}\left(u_0(0,y)y\right) &=& \frac{\partial u_0}{\partial y}(0,y)y +u_0(0,y) \;=\; 0.
\end{eqnarray}
This condition has also to be satisfied for the solution $y=y_k$ of (\ref{zeroimaxcond1}) which make this equation to the identity (\ref{yk0}) that means that
\begin{eqnarray}
\frac{\partial u_0}{\partial y}(0,y_k)y_k +u_0(0,y_k) &=& \frac{\partial u_0}{\partial y}(0,y_k)y_k +\frac{n\pi}{y_k} \;=\; 0,
\end{eqnarray}
has to be identically satisfied. Moreover, if we apply the operator $\frac{\sin^m\left(x\frac{\partial}{\partial y}\right)}{x^{m+1}}$ to condition (\ref{zeroscond1}) we obtain
\begin{eqnarray}
\frac{\sin^{m+1}\left(x\frac{\partial}{\partial y}\right)}{x^{m+1}}\left(u_0(0,y)y\right) &=& \left(\frac{\sin\left(x\frac{\partial}{\partial y}\right)}{x\frac{\partial}{\partial y}}\right)^{m+1} \frac{\partial^{m+1}}{\partial y^{m+1}}\left(u_0(0,y)y\right) \;=\; 0,\;\; (m=0,1,2,\ldots),\qquad\; \label{zeroscond4}
\end{eqnarray}
and by multiplication with the inverse operator to the nonsingular operator $\left(\frac{x\frac{\partial}{\partial y}}{\sin\left(x\frac{\partial}{\partial y}\right)}
\right)^{m+1} $ we find
\begin{eqnarray}
\frac{\partial^{m+1}}{\partial y^{m+1}}\left(u_0(0,y)y\right) \;=\; 0,\quad (m=0,1,2,\ldots).\label{condzerosnplus1}
\end{eqnarray}
All these conditions have to be satisfied for the solutions $y=y_k$ with $x\neq 0$ in (\ref{zeroimaxcond2}), i.e.
\begin{eqnarray}
\frac{\partial^{m+1}}{\partial y^{m+1}}\left\{\left(u_0(0,y)y\right)\right\}_{y=y_k} \;=\; 0,\quad (m=0,1,2,\ldots). \label{dnyk0}
\end{eqnarray}
The same conditions follow also from (\ref{zeroscond2}) combined with (\ref{zeroscond2}) by Taylor series expansion with respect to $x\frac{\partial}{\partial y}$ for $x\neq 0$ according to
\begin{eqnarray}
n\pi &=&\exp\left(\pm \I x\frac{\partial}{\partial y}\right)\left(u_0(0,y)y\right) \;=\; u_0(0,y)y+\sum_{m=0}^\infty\frac{\left(\pm \I x\right)^{m+1}}{(m+1)!}\frac{\partial^{m+1}}{\partial y^{m+1}}\left(u_0(0,y)y\right),
\end{eqnarray}
by setting all sum terms proportional to $x^{m+1}$ equal to zero.

If we make now a Taylor series expansion of the function $u_0(0,y)y$ in the neighborhood $y=y_k$ of a solution which obeys all conditions (\ref{yk0}) and (\ref{dnyk0}) then we find
\begin{eqnarray}
u_0(0,y)y &=& u_0(0,y_k)y_k +\sum_{m=0}^{\infty}\frac{1}{(m+1)!}\frac{\partial^{m+1}}{\partial y^{m+1}}\left\{\left(u_0(0,y)y\right)\right\}_{y=y_k}\left(y-y_k\right)^{m+1} \;=\; n\pi.
\end{eqnarray}
Thus we can find zeros for $x\neq 0$ that means off the imaginary axis if the mean value function $u_0(0,y)$ on the imaginary axis possess one of the forms
\begin{eqnarray}
u_0(0,y)y = n\pi, \quad \Leftrightarrow \quad u_0(0,y) = \frac{n\pi}{y},\quad (n=0,\pm 1,\pm 2,\ldots),
\end{eqnarray}
for a certain integer $n$. According to (\ref{uvxycossin}) the whole mean-value functions $u_0(x,y)$ and $v_0(x,y)$ are then
\begin{eqnarray}
u_0(x,y) &=& \cos\left(x\frac{\partial}{\partial y}\right)\frac{n\pi}{y} \;=\; \frac{n\pi}{2}\left\{\exp\left(\I x\frac{\partial}{\partial y}\right) +\exp\left(-\I x\frac{\partial}{\partial y}\right)\right\}\frac{1}{y} \nonumber\\ &=& \frac{n\pi}{2}\left(\frac{1}{y+\I x}+\frac{1}{y-\I x} \right) \;=\; n\pi \frac{y}{x^2+y^2}, \nonumber\\ v_0(x,y) &=& -\sin\left(x\frac{\partial}{\partial y}\right)\frac{k\pi}{y} \;=\; -\frac{n\pi}{\I 2}\left\{\exp\left(\I x\frac{\partial}{\partial y}\right) -\exp\left(-\I x\frac{\partial}{\partial y}\right)\right\}\frac{1}{y} \nonumber\\ &=& -\frac{n\pi}{\I 2}\left(\frac{1}{y+\I x}-\frac{1}{y-\I x} \right) \;=\; n\pi \frac{x}{x^2+y^2},
\end{eqnarray}
or in compact form
\begin{eqnarray}
w_0(z) &=& u_0(x,y)+\I v_0(x,y) \;=\; n\pi\frac{y+\I x}{x^2+y^2} \;=\; \I\frac{n\pi}{x+\I y},\quad \Leftrightarrow \quad w_0(z)z \;=\; \I n\pi.
\end{eqnarray}
If we insert $w_0(z)z= \I n\pi$ into equation (\ref{secmeanvth2compl}) then we get ${\mit \Xi}(z) = 0$ for all $n=0,\pm 1,\pm2,\ldots$. This means that all conditions for zeros with $x\neq 0$ together do not lead to a solution for arbitrary ${\mit \Xi}(z) \neq 0$. Thus we have proved that all zeros of Xi functions ${\mit \Xi}(z)$ lie on the imaginary axis $x=0$.

Recognizing equations (\ref{zeroscond1}) and (\ref{zeroscond2}) as correct ones one may find modifications of the given proof of impossibility of solutions for zeros with $x\neq 0$. For example, one can make Taylor series expansions and obtain from equation (\ref{zeroscond2})
\begin{eqnarray}
u_0(0,y)y+\sum_{m=1}^{\infty}\frac{(-1)^m x^{2m}}{(2m)!}\frac{\partial^{2m}}{\partial y^{2m}}\left(u_0(0,y)y\right) &=& n\pi, \quad (n=0,\pm 1,\pm 2,\ldots),\label{taylu00y1}
\end{eqnarray}
and from (\ref{zeroscond1})
\begin{eqnarray}
\sum_{m=0}^{\infty}\frac{(-1)^m x^{2m+1}}{(2m+1)!}\frac{\partial^{2m+1}}{\partial y^{2m+1}}\left(u_0(0,y)y\right) &=& 0. \label{taylu00y2}
\end{eqnarray}
The right-hand sides of these equations are independent of variable $x$ and therefore the left-hand sides must it too. On the imaginary axis for $x=0$ this leads to the condition (\ref{zeroimaxcond2}) as the only condition which determines the zeros on this axis. For $x\neq 0$
one has as additional condition the vanishing of the coefficients in front of powers of $x^{n+1},\,(n=0,1,2,\ldots)$ in (\ref{taylu00y1}) and (\ref{taylu00y2}) that taken together leads to the conditions (\ref{condzerosnplus1}). The further discussion of the impossibility of this case is the same as before.

A third variant to show the impossibility for zeros in case of $x\neq 0$ is to make the transition to the Fourier transform of $u_0(0,y)y$
and to solve the equation arising from \ref{zeroscond3}) by generalized functions and then making the inverse Fourier transformation. One may show then that this is not compatible with the general solution of (\ref{zeroscond2}) which determines the position of the zeros. We do not present this variant here.

We have now finally proved that all Xi functions ${\mit \Xi}(z)$ of the form (\ref{basiccapxi0}) for which the second mean-value theorem is applicable (function ${\mit \Omega}(u)$ positively semi-definite and non-increasing (or also negatively semi-definite and non-decreasing)
may possess zeros only on the imaginary axis.

\setcounter{chapter}{6}
\setcounter{equation}{0}
\section*{6. Consequences for proof of the Riemann hypothesis}

The given proof for zeros only on the imaginary axis $x=0$ for the considered Xi function ${\mit \Xi}(z)={\mit \Xi}(x+\I y)$ includes as special case the function ${\mit \Omega}(u)$ to the Riemann hypothesis which is given in (\ref{defcapomega}). However, it includes also the whole class of modified Bessel functions of imaginary argument ${\rm I}_{\nu}(z)$ which possess zeros only on the imaginary axis and if we make the substitution $z\leftrightarrow \I z$ also the usual Bessel function ${\rm J}_{\nu}(z)$ which possess zeros only on the real axis.

We may ask about possible degeneracies of the zeros of the Xi functions ${\mit \Xi}(z)$ on the imaginary axis $z=\I y$. Our proof does not give a recipe to see whether such degeneracies are possible or not. In case of the Riemann zeta function ${\mit \Xi}(z) \leftrightarrow \xi(s) \leftrightarrow \zeta(s)$ one cannot expect a degeneracy because the countable number of all nontrivial zeros are (likely)
irrational (transcendental, proof?) numbers but we do not know a proof for this.

For ${\mit \Xi}(z)$ as an entire function one may pose the question of its factorization with factors of the form $1-\frac{z}{z_n} $ where $z_n$ goes through all roots where in case of degeneracy the same factors are taken multiple times according to the degeneracy.
It is well known that an entire function using its ordered zeros $z_n,\, (|z_n|\le |z_{n+1}|)$ can be represented in Weierstrass product form multiplied by an exponential function $\E^{h(z)}$ with an entire function function $h(z)$ in the exponent with the result that $\E^{h(z)}$ is an entire function without zeros. This possesses the form (e.g., \cite{lang})
\begin{eqnarray}
{\mit \Xi}(z) &=& \E^{h(z)} z^{m}\prod_{n}\left(1-\frac{z}{z_{n}}\right)\exp\left(P_{k_n}\left(\frac{z}{z_n}\right)\right),\label{wei}
\end{eqnarray}
with a polynomial $P_{k}(w)$ of degree $k$ which depending on the roots $z_n$ must be appropriately chosen to guarantee the convergence of the product. This polynomial is defined by first $k$ sum terms in the Taylor series for $-\log(1-w)$\footnote{Sometimes our $P_{k-1}(w)$ is denoted by $P_{k}(w)$.}
\begin{eqnarray}
P_k(w) & \equiv & w+\frac{w^2}{2}+\frac{w^3}{3}+\ldots+\frac{w^{k}}{k} \;=\; -\log(1-w)-\sum_{l=k+1}^{\infty}\frac{w^l}{l}.
\end{eqnarray}
By means of these polynomials the Weierstrass factors are defined as the functions
\begin{eqnarray}
E_k(w) &\equiv & (1-w)\exp\left(P_k(w)\right),
\end{eqnarray}
from which follows
\begin{eqnarray}
\log\left(E_k(w)\right) \;=\; \log(1-w)-\left(\log(1-w)+\sum_{l=k+1}^{\infty}\frac{w^l}{l}\right) \;=\; -\sum_{l=k+1}^{\infty}\frac{w^l}{l}.
\end{eqnarray}
From this form it is seen that $E_k(w)$ possesses the following initial terms of the Taylor series
\begin{eqnarray}
E_k(w) &=& (1-w)\exp\left(-\sum_{l=k+1}^{\infty}\frac{w^{l}}{l}\right) \;=\; 1-\frac{w^{k+1}}{k+1} +\ldots\;,
\end{eqnarray}
and is a function with a zero at $w =1$ but with a Taylor series expansion which begins with the
terms $1-\frac{w^{k+1}}{k+1}$.
Hadamard made a precision of the Weierstrass product form by connecting the degree $k_n$ of the polynomials in (\ref{wei}) with the order $\rho$ of growth of the entire function and showed that $k_n$ can be chosen independently of the $n$-th root $z_n$ by $k_n \rightarrow k \ge \rho -1$
The order of ${\mit \Xi}(z)$ which is equal to $1$ is not a strict order $\rho = 1$ (for this last notion see \cite{lang}). However, this does not play a role in the Hadamard product representation of ${\mit \Xi}(z)$ and the polynomials
$P_{k_n}(w)$ in (\ref{wei}) can be chosen as $P_0(w)$ that means equal to $0$ according to $k_n= k =\rho -1$. The entire function $h(z)$ in the exponent in (\ref{wei}) can be only a constant since in other case it would introduce a higher growth of ${\mit \Xi}(z)$. Thus the product representation of ${\mit \Xi}(z)$ possesses the form
\begin{eqnarray}
{\mit \Xi}(z) &=& {\mit \Xi}(0){{\prod}'}_{n=-\infty}^{+\infty} \left(1-\frac{z}{z_n}\right) \;=\; \int_{0}^{+\infty}du\,{\mit \Omega}(u) \prod_{n=1}^{+\infty}
\left(1-\frac{z^2}{z_n^2}\right) \nonumber\\ &=& {\mit \Omega}_{0} \prod_{n=1}^{+\infty}\left(1+\frac{z^2}{y_n^2}\right), \quad {\mit \Xi}(0) \;=\; {\mit \Omega}_0 \;\equiv \; \int_{0}^{+\infty}du\,{\mit \Omega}(u) \;=\; \xi\left(\frac{1}{2}\right) \;=\; 0.49712...\;,\label{prodxi}
\end{eqnarray}
where we took into account the symmetry $z_{-n} = -z_n =(z_n)^*$ of the zeros and the proof $z_n=\I y_n$ that all zeros lie on the imaginary axis and a zero $z_0=0$ is absent. With ${\mit \Omega}_0$ we denoted the first moment of the function ${\mit \Omega}(u)$.

Formula (\ref{prodxi}) in connection with his hypothesis was already used by Riemann in \cite{rie}
and later proved by von Mangoldt where the product representation of entire functions by Weierstrass which was later stated more precisely by Hadamard plays a role. There is another formula for an approximation to the number of nontrivial zeros of $\zeta(s)$ or $\xi(s)$ which in application to the number of zeros $N(Y)$ of ${\mit \Xi}(z)$ on the imaginary axis $z=\I y$ (critical line) in the interval between $y=0$ and $y=Y$. It takes on the form ($Y$ for ${\mit \Xi}(z)$ is equivalent to usual $T$ for $\zeta(s)$)
\begin{eqnarray}
N(Y) &=& \int_{0}^{Y}dy\,\nu(y) \;\approx\;\frac{Y}{2\pi}\log\left(\frac{Y}{2\pi}\right)-\frac{Y}{2\pi},\quad (Y \gg 0),\label{xinumbzeros}
\end{eqnarray}
with the logarithmically growing density
\begin{eqnarray}
\nu(y) &\approx & \frac{1}{2\pi}\log\left(\frac{y}{2\pi}\right),\quad (y \gg 1). \label{xidenszeros}
\end{eqnarray}
As long as the Riemann hypothesis was not proved it was formulated for the critical strip $ 0 \le \sigma \le 1$ of the complex coordinate $s=\sigma +\I t$ in $\xi(s)$ parallel to the imaginary axis and with $t$ between $t=0$ and $t=T$ (with $T$ equal to our $Y$ in (\ref{xinumbzeros})). It was already suggested by Riemann \cite{rie} but not proved in detail there and was later proved by von Mangoldt in 1905. A detailed proof by means of the argument principle can be found in \cite{borw}. The result of Hardy (1914) (cited in \cite{edw}) that there exist an infinite number of zeros on the critical line is a step to the full proof of the Riemann hypothesis. Section 4 of present article may be considered as involving such proof of this last statement.

We have now proved that functions ${\mit \Xi}(z)$ defined by integrals of the form (\ref{basiccapxi0}) with non-increasing functions ${\mit \Omega}(u)$ which decrease in infinity sufficiently rapidly in a way that ${\mit \Xi}(z)$ becomes an entire function of $z$ possess zeros only on the imaginary axis $z=\I y$. This did not provide a recipe to see in which cases all zeros on the imaginary axis are simple zeros but it is unlikely that within a countable sequence of irregularly chosen (probably) transcendental numbers (the zeros) two of them are coincident (it seems to be difficult to formulate last statement in a more rigorous way). It also did not provide a direct formula for the number of zeros in an interval $\left[(0,0),(0,Y)\right]$ from zero to $Y$ on the imaginary axis or of its density there but, as already said, Riemann \cite{rie} suggested for this an approximate formula and von Mangoldt proved it

The proof of the Riemann hypothesis is included as the special case (\ref{defcapomega}) of the function ${\mit \Omega}(u)$ into a wider class of functions with an integral representation of the form (\ref{basiccapxi0}) which under the discussed necessary conditions allowing the application of the second mean-value theorem of calculus possess zeros only on the imaginary axis.
The equivalent forms (\ref{capxirep6}) and (\ref{capxirep7}) of the integral (\ref{basiccapxi0}) where the functions, for example ${\mit \Omega}^{(1)}(u)$, are no more generally non-increasing suggest that conditions for zeros only on the imaginary axis are existent for more general cases than such prescribed by the second mean-value theorem. A certain difference may happen, for example, for $z=0$ because powers of it are in the denominators in the representations in (\ref{capxirep7}).

\setcounter{chapter}{7}
\setcounter{equation}{0}
\section*{7. Graphical illustration of mean-value parameters to Xi function for the Riemann hypothesis}

\begin{figure}[h]
\includegraphics[width=16.0cm]{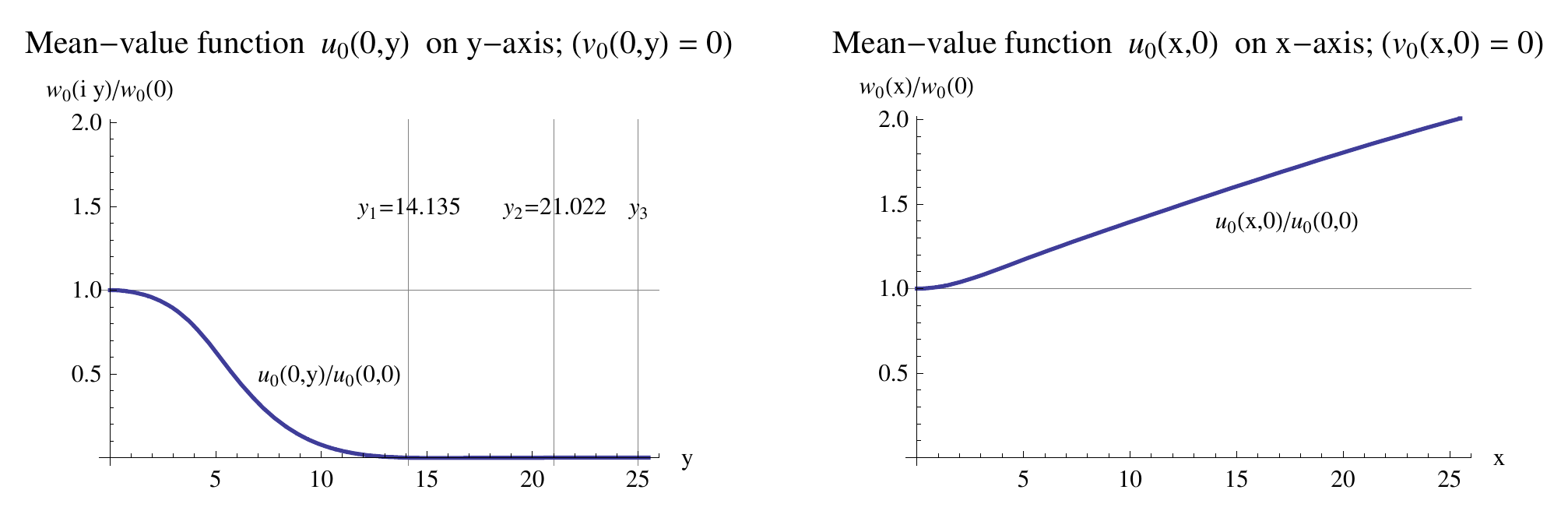}\hspace{1cm}
\caption{\foot Mean value parameters $w_0(\I y) = u_0(0,y)$ and $w_0(x) = u_0(x,0)$ for the Xi function in the proof of the Riemann hypothesis. \newline \scriptsize{It is not to see in the chosen scale that the curve $u_0(0,y)$ goes beyond the $y$-axis and oscillates around it due to extremely rapid vanishing of the envelope of $u_0(0,y)$ with increasing $y$ but we do not resolves this here by additional graphics because this behavior is better to see in the later considered case of modified Bessel functions. Using (\ref{w0zriem}) we calculate numerically $w_0(0)=u_0(0,0) \approx 0.27822$ that is the value which we call the optimal value for the moment series expansion (see Section 11). The part in the second partial figure which at the first glance looks like a straight line as asymptote is not such.}}
\end{figure}

\begin{figure}[h]
\includegraphics[width=16.0cm]{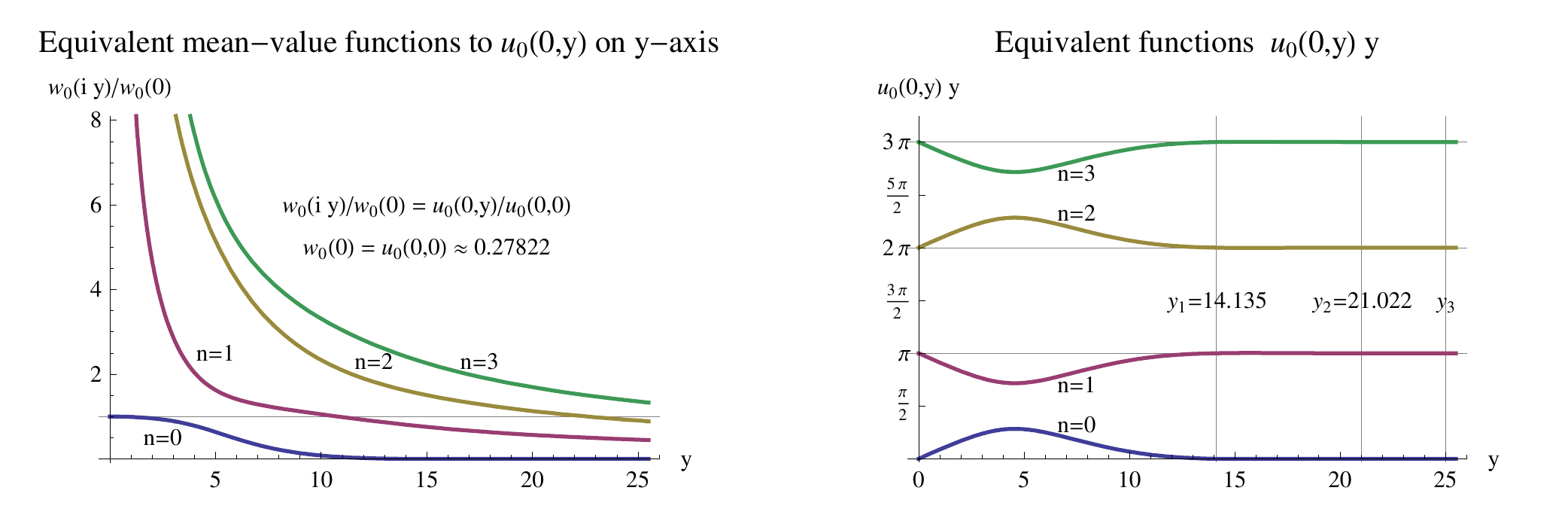}\hspace{1cm}
\caption{\foot Mean value parameters in the proof of the Riemann hypothesis. \newline \scriptsize{On the left-hand side there are shown the mean value parameters $\frac{w_0(\I y)}{w(0)} = \frac{u_0(0,y)}{u_0(0,0)}$ for the Xi function to the Riemann hypothesis if we do not take the values of the function ${\rm \arcsin(t)}$ in the basic range $-\frac{\pi}{2} \le t \le \frac{\pi}{2}$ but in equivalent ranges according to (\ref{arcsinequ}). On the right-hand side are
shown the corresponding functions $u_0(0,y)y$ which according to $\cos(x+n\pi)=(-1)^n\sin(x)$ and the condition for zeros $\cos\left(u_0(0,y)y\right)=0$ lead to equivalent ranges
$k\pi =u_0(0,y)y \cong u_0(0,y)y +n\pi= (k+n)\pi, (k=\pm 1, \pm 2,\ldots;\,n=0,\pm 1,\pm2,\ldots)$ (see (\ref{zerosimax})) determine the zeros of the Xi function on the imaginary axis. We see that the multi-valuedness of the ${\rm arcsin(x)}$ function does not spoil a unique result for the zeros because every branch find the corresponding $n$ of $n\pi$ where then all zeros lie. Due to extremely rapid decrease of the function $u_0(0,y)$ with increasing $y$ this is difficult to see (position of first three zero at $y_1\approx 14.1,\,y_2 \approx 21.0,\,y_3 \approx 25.0$ is shown) but if we separate small intervals of $y$ and enlarge the range of values for $y_0(0,y)$ this becomes visible (similar as in Fig. 2.3). We do not make this here because this effect is better visible for the modified Bessel functions which we intend to consider at another place.}}
\end{figure}

To get an imagination how the mean-value function $w_0(z) =w_0(x+\I y)$ looks like we calculate it for the imaginary axis and for the real axis for the case of the function ${\mit \Omega}(u)$ in (\ref{defcapomega}) that is possible numerically. From the two equations for general $z$ and for $z=0$
\begin{eqnarray}
{\mit \Omega}(0) \frac{\,{\rm sh}\left(w_0(z)z\right)}{z} &=& \int_{0}^{+\infty}du\,{\mit \Omega}(u)\,{\rm ch}(uz), \nonumber\\ {\mit \Omega}(0)w_0(0) &=& \int_{0}^{+\infty}du\,{\mit \Omega}(u),
\end{eqnarray}
follows
\begin{eqnarray}
w_0(z) &=& \frac{1}{z}\,{\rm Arsh}\left(\frac{z}{{\mit \Omega}(0)}\int_{0}^{+\infty}du\,{\mit \Omega}(u)\,{\rm ch}(uz)\right), \nonumber\\ w_{0}(0) &=& \frac{1}{{\mit \Omega}(0)}\int_{0}^{+\infty}du\,{\mit \Omega}(u),  \label{w0zriem}
\end{eqnarray}
with the two initial terms of the Taylor series
\begin{eqnarray}
w_0(z) &=& \frac{1}{z}\left\{\frac{z}{{\mit \Omega}(0)}\int_{0}^{+\infty}du\,{\mit \Omega}(u)\,{\rm ch}(uz) -\frac{1}{6}\left(\frac{z}{{\mit \Omega}(0)}\int_{0}^{+\infty}du\,{\mit \Omega}(u)\,{\rm ch}(uz)\right)^3 +\ldots\right\},\label{w0zriemtaylor}
\end{eqnarray}
and with the two initial terms of the asymptotic series
\begin{eqnarray}
w_0(z) \;=\; \frac{1}{z}\left\{\log\left(2\frac{z}{{\mit \Omega}(0)}\int_{0}^{+\infty}du\,{\mit \Omega}(u)\,{\rm ch}(uz)\right) +\frac{1}{4}\left(\frac{z}{{\mit \Omega}(0)}\int_{0}^{+\infty}du\,{\mit \Omega}(u)\,{\rm ch}(uz)\right)^{-2}+\ldots \right\}. \label{w0zriemasympt}
\end{eqnarray}

From (\ref{w0zriem}) follows
\begin{eqnarray}
\frac{w_0(z)}{w_0(0)} &=& \frac{{\mit \Omega}(0)}{z\int_{0}^{+\infty}du\,{\mit \Omega}(u)}\,{\rm Arsh}\left(\frac{z}{{\mit \Omega}(0)}\int_{0}^{+\infty}du\,{\mit \Omega}(u)\,{\rm ch}(uz)\right).\label{w0zriem1}
\end{eqnarray}
This can be numerically calculated from the explicit form (\ref{defcapomega}) of ${\mit \Omega}(u)$. For $y=0$ and for $x=0$  (and only for these cases) the function $\frac{w_0(z)}{w_0(0)}$ is real-valued, in particular, for $y=0$
\begin{eqnarray}
\frac{w_0(x)}{w_0(0)} &=& \frac{{\mit \Omega}(0)}{x\,\int_{0}^{+\infty}du\,{\mit \Omega}(u)}\,{\rm Arsh}\left(\frac{x}{{\mit \Omega}(0)}\int_{0}^{+\infty}du\,{\mit \Omega}(u)\,{\rm ch}(ux)\right) \nonumber\\ &=& \frac{\int_{0}^{+\infty}du\,{\mit \Omega}(u)\,{\rm ch}(ux)}{\int_{0}^{+\infty}du\,{\mit \Omega}(u)}\left\{1-\frac{1}{6}\left(\frac{x}{{\mit \Omega}(u)}\int_{0}^{+\infty}du\,{\mit \Omega}(u)\,{\rm ch}(ux)\right)^2+\ldots \right\} \nonumber\\ &=& \frac{u_0(x,0)}{u_0(0,0)},\label{meanvalwxriem}
\end{eqnarray}
and for $x=0$
\begin{eqnarray}
\frac{w_0(\I y)}{w_0(0)} &=& \frac{{\mit \Omega}(0)}{y\,\int_{0}^{+\infty}du\,{\mit \Omega}(u)}\,{\rm arcsin}\left(\frac{y}{{\mit \Omega}(0)}\int_{0}^{+\infty}du\,{\mit \Omega}(u)\,{\rm cos}(uy)\right) \nonumber\\ &=&
\frac{\int_{0}^{+\infty}du\,{\mit \Omega}(u)\,{\rm cos}(uy)}{\int_{0}^{+\infty}du\,{\mit \Omega}(u)}\left\{1+\frac{1}{6}\left(\frac{y}{{\mit \Omega}(u)}\int_{0}^{+\infty}du\,{\mit \Omega}(u)\,{\rm cos}(uy)\right)^2+\ldots \right\} \nonumber\\ &=& \frac{u_0(0,y)}{u_0(0,0)} ,\qquad\label{meanvalwyriem}
\end{eqnarray}
where we applied the first two terms of the Taylor series expansion of ${\rm arcsin}(x)$ in powers of $x$. A small problem is here that we get the value for this multi-valued function in the range $-\frac{\pi}{2} \le {\rm arcsin}(x) \le +\frac{\pi}{2}$. Since $\frac{1}{x}\, {\rm arcsin}(x)$ is
an even function with only positive coefficients in its Taylor series the term in braces is in every case positive that becomes important below.

The two curves which we get for $\frac{w_0(x)}{w_0(0)}$ and for $\frac{w_0(\I y)}{w_0(0)}$ are shown in Fig. 7.4. The function for $w_0(x,0)$ on the real axis $y=0$ (second partial picture) is not very exciting. The necessary condition $xu_0(x,0)=0$ (see (\ref{neccondzeros})) can be satisfied only for $x=0$ but it is easily to see from ${\mit \Xi}(0) = \int_{0}^{+\infty}du\;{\mit \Omega}(u) \approx 1.7868 \neq 0$ that there is no zero. For the function $w_0(0,y)$ on the imaginary axis $x=0$ the necessary condition $yv_0(0,y)=0$ (see (\ref{neccondzeros})) is trivially satisfied since $v_0(0,y)=0$ and does not restrict the solutions for zeros. In this case only the sufficient condition $yu_0(0,y) = n\pi,(n=0,\pm 1,\pm2,\ldots)$ determines the position of the zeros on the imaginary axis.  The first two pairs of zeros are at $y_1\approx \pm 14.135,\;y_2 \approx \pm 21.022$ and the reason that we do not see them in Fig. 7.4 is the rapid decrease of the function $u_0(0,y)$ with increasing $y$. If we enlarge this range we see that the curve goes beyond the $y$-axis after the first root at $14.135$ of the Xi function. As a surprise for the second mean-value method we see that the parameter $u_0(0,y)$ becomes oscillating around this axis.
This means that the roots which are generally determined by the equation $yu_0(0,y)=n\pi$ (see (\ref{zerosequ})) are determined here by the value $n=0$ alone. The reason for this is the multi-valuedness of the ArcSine function according to
\begin{eqnarray}
{\rm arcsin}(x) & \cong & n\pi+(-1)^n {\rm arcsin}(x),\quad \left(-1 \le x \le +1,\quad -\frac{\pi}{2} \le {\rm arcsin}(x) \le +\frac{\pi}{2},\quad n \in {\mathbb Z} \right).\qquad
\label{arcsinequ}
\end{eqnarray}
If we choose the values for the ${\rm arcsin(x)}$-function not in the basic interval $-\frac{\pi}{2} \le x \le +\frac{\pi}{2}$ for which the Taylor series provides the values but from other equivalent intervals according to (\ref{arcsinequ}) we get other curves for $u_0(0,y)$ and $yu_0(0,y)$ from which we also may determine the zeros (see Fig. 7.5), however, with other values $n$ in the relation $yu_0(0,y)=n\pi,\,(n=0,\pm 1,\pm 2,\ldots)$ and the results are invariant with respect to the multi-valuedness. This is better to see in case of the modified Bessel functions for which the curves vanish less rapidly with increasing $y$ as we intend to show at another place. All these considerations do not touch the proof of the non-existence of roots off the imaginary axis but should serve only for better understanding of the involved functions. It seems that the specific phenomenons of the second mean-value theorem (\ref{secmeanvth}) if the functions $g(u)$ there are oscillating functions (remind, only continuity is required) are not yet well illustrated in detail.

We now derive a few general properties of the function $u_0(0,y)$ which can be seen in the Figures. From (\ref{realsecmeanvth}) written in the form and by Taylor series expansion according to
\begin{eqnarray}
\int_{0}^{+\infty}du\,{\mit \Omega}(u)\cos(u y) &=& {\mit \Omega}(0)u_0(0,y) \frac{\sin\left(u_0(0,y)y\right)}{u_0(0,y)y} \nonumber\\ &=& {\mit \Omega}(0)u_0(0,y) \left(1-\frac{y^2}{6}\left(u_0(0,y)\right)^2 +\ldots \right),\label{secmeanvth2real2}
\end{eqnarray}
follows from the even symmetry of the left-hand side that $u_0(0,y)$ also has to be a function of the variable $y$ with even symmetry (notation $\frac{\dis{\partial^n u_0(0,y)}}{\dis{\partial y^n}} \equiv u_0^{(n)}(0,y)$)
\begin{eqnarray}
u_0(0,y) &=& +u_0(0,-y),\label{secmeanvth2real3}
\end{eqnarray}
with the consequence
\begin{eqnarray}
u_0(0,y) \;=\; u_0(0,0)+\sum_{m=1}^{\infty}\frac{u_0^{(2m)}(0,0)}{(2m)!}y^{2m},\quad u_0^{(2m+1)}(0,0)\;=\; 0.\qquad \label{secmeanvth2real30}
\end{eqnarray}
Concretely, we obtain by $n$-fold differentiation of both sides of (\ref{secmeanvth2real2}) at $y=0$ for the first coefficients of the Taylor series
\begin{eqnarray}
\int_{0}^{+\infty}du\,{\mit \Omega}(u) &=& {\mit \Omega}(0)u_0(0,0), \quad 0 \;=\; {\mit \Omega}(0)u_0^{(1)}(0,0), \nonumber\\ -\int_{0}^{+\infty}du\,{\mit \Omega}(u)u^2 &=& {\mit \Omega}(0)\left(u_0^{(2)}(0,0)-\frac{1}{3}\left(u_0(0,0)\right)^3\right),
\end{eqnarray}
from which follows
\begin{eqnarray}
u_0(0,0) &=& \frac{1}{{\mit \Omega}(0)}\int_{0}^{+\infty}du\,{\mit \Omega}(u), \quad u_0^{(1)}(0,0) \;=\; 0, \quad  \nonumber\\ u_0^{(2)}(0,0) &=& \underbrace{-\frac{1}{{\mit \Omega}(0)}\int_{0}^{+\infty}du\,{\mit \Omega}(u)u^2}_{<\,0} +\underbrace{\frac{1}{3}\left(\frac{1}{{\mit \Omega}(0)}\int_{0}^{+\infty}du\,{\mit \Omega}(u)\right)^3}_{>\,0}.
\end{eqnarray}
Since the first sum term on the right-hand side is negative and the second is positive it depends from their values whether or not $u_0^{(2)}(0,0)$ possesses a positive or negative value. For the special function ${\mit \Omega}(u)$ in (\ref{defcapomega}) which plays a role in the Riemann hypothesis we find approximately
\begin{eqnarray}
&& {\mit \Omega}_0 \;\equiv\; \int_{0}^{+\infty}du\,{\mit \Omega}(u) \;\approx \; 0.497121,\quad 2!{\mit \Omega}_2 \;\equiv \; \int_{0}^{+\infty}du\,{\mit \Omega}(u)u^2 \;\approx \; 0.0229719,\nonumber\\ && {\mit \Omega}(0) \;\approx\; 1.78679,\quad
u_0(0,0) \;\approx \; 0.278220,\quad u_0^{(2)}(0,0) \;\approx\; -0.00567784,\label{momentnumer}
\end{eqnarray}
meaning that the second coefficient in the expansion of $u_0(0,y)$ in a Taylor series in powers of $y$ is negative that can be seen in the first part of Fig. 7.4. However, as we have seen the proof of the Riemann hypothesis is by no means  critically connected with some numerical values.

In principle, the proof of the Riemann hypothesis is accomplished now and illustrated and we will stop here. However, for a deeper understanding of the proof it is favorable to consider some aspects of the proof such as, for example, analogues to other functions with a representation of the form (\ref{basiccapxi0}) and with zeros only on the imaginary axis and some other approaches although they did not lead to the full proof that,
however, we cannot make here.

\setcounter{chapter}{8}
\setcounter{equation}{0}
\section*{8. Equivalent formulations of the main theorems in a summary}

In present article we proved the following main result\\[1mm]
{\large{\bf Theorem 1:}}\\
Let ${\mit \Omega}(u)$ be a real-valued function of variable $u$ in the interval $0 \le u < +\infty$ which is positive semi-definite in this interval and non-increasing and is rapidly vanishing in infinity, more rapidly than any exponential function $\exp\left(-\lambda u\right)$, that means
\begin{eqnarray}
{\mit \Omega}(u) \; \ge 0, \quad {\mit \Omega}^{(1)}(u) &\le & 0,\quad (0\le u <+\infty),\quad \lim_{u\rightarrow +\infty}\frac{{\mit \Omega}(u)}{\exp\left(-\lambda u\right)} \;= \;0,\quad (\lambda \in {\mathbb R}^{+}).
\end{eqnarray}
Then the following integral with arbitrary complex parameter $z=x+\I y$
\begin{eqnarray}
{\mit \Xi}(z) &=& \int_{0}^{+\infty}du\,{\mit \Omega}(u) \,{\rm ch}(uz), \quad (z \in {\mathbb C}), \label{basiccapxi2}
\end{eqnarray}
is an entire function of $z$ with possible zeros $z=z_0=x_0+\I y_0$ only on the imaginary axis $x=0$ that means
\begin{eqnarray}
{\mit \Xi}(z_0) &=& 0, \quad \Rightarrow \quad z_0\;=\; \I y_0,\quad ({\rm or}\;\; x_0=0).
\end{eqnarray}
{\large{\bf Proof:}}\\[1mm]
The proof of this theorem for non-increasing functions ${\mit \Omega}(u)$ takes on Sections 3--5 of this article. The function ${\mit \Omega}(u)$ in (\ref{defcapomega}) satisfies these conditions and thus provides a proof of the Riemann hypothesis.
\\[1mm]
{\large{\bf Remark:}}\\
A analogous theorem is obviously true by substituting in (\ref{basiccapxi2}) ${\rm ch}(u) \leftrightarrow \cos(u)$ and by interchanging the role of the imaginary and of the real axis $y\leftrightarrow x$.
Furthermore, a similar theorem with a few peculiarities is true for substituting ${\rm ch}(uz)$ in (\ref{basiccapxi2}) by ${\rm sh}(uz)$ .\\
\indent

{\bf{Theorem 1}} can be formulated in some equivalent ways which lead to interesting consequences\footnote{Some of these equivalences now formulated as consequences originate from trials to prove the Riemann hypothesis in other way.}. The Mellin transformation ${\hat{f}}(s)$ of an arbitrary function $f(t)$ together with its inversion is defined by \cite{erd,zem,brpr,bert}
\begin{eqnarray}
{\hat{f}}(s) \;=\; \int_{0}^{+\infty}dt\,f(t)t^{s-1},\quad f(t)\;=\; \frac{1}{\I 2\pi}\int_{c-\I \infty}^{c+\I \infty}ds\,{\hat{f}}(s)t^{-s},\label{mellintransform}
\end{eqnarray}
where the real value $c$ has only to lie in the convergence strip for the definition of ${\hat{f}}(s)$ by the integral.
Formula (\ref{basiccapxi2}) is an integral transform of the function ${\rm ch}(z)$ and can be considered as the application of an integral operator to the function ${\rm ch}(z)$ which using the Mellin transform ${\hat{\mit \Omega}}(s)$ of the function ${\mit \Omega}(u)$ can be written in the following convenient form
\begin{eqnarray}
{\mit \Xi}(z) &=& {\hat{\mit \Omega}}\left(\frac{\partial}{\partial z}z\right) \,{\rm ch}(z).\label{theoremmellin}
\end{eqnarray}
This is due to
\begin{eqnarray}
\int_{0}^{+\infty}du\,{\mit \Omega}(u) \,{\rm ch}(uz) &=& \int_{0}^{+\infty}du\,{\mit \Omega}(u) u^{z\frac{\partial}{\partial z}}\,{\rm ch}(z) \;=\; \int_{0}^{+\infty}du\,{\mit \Omega}(u) u^{\frac{\partial}{\partial z}z-1}\,{\rm ch}(z),
\end{eqnarray}
where $u^{z\frac{\partial}{\partial z}}$ is the operator of multiplication of the argument of an arbitrary function $g(z)$ by the number $u$, i.e. it transforms as follows
\begin{eqnarray}
g(z) \;\rightarrow \; g(uz) &=& u^{z\frac{\partial}{\partial z}} g(z),
\end{eqnarray}
according to the following chain of conclusions starting from the property that all functions $z^n,\,(n=0,1,2,\ldots)$ are eigenfunctions of $z\frac{\partial}{\partial z}$ to eigenvalue $n$
\begin{eqnarray}
&& z\frac{\partial}{\partial z}z^{n} \;=\; nz^{n},\quad \Rightarrow \quad f\left(z\frac{\partial}{\partial z}\right) z^n \;=\; f(n)z^n, \quad \Rightarrow \quad \exp\left(\lambda z\frac{\partial}{\partial z}\right)z^n \;=\; \left(\E^{\lambda} z\right)^n, \quad \Rightarrow \nonumber\\ && \exp\left(\lambda z\frac{\partial}{\partial z}\right)g(z) \;=\; g\left(\E^\lambda z\right),\quad \Rightarrow \quad u^{z\frac{\partial}{\partial z}} g(z) \;=\; g(uz),\quad (\E^{\lambda}\equiv u).
\end{eqnarray}
This chain is almost obvious and does not need more explanations. The operators $u^{z\frac{\partial}{\partial z}}$ are linear operators in linear spaces depending on the considered set of numbers $u$.

Expressed by real variables $(x,y)$ and by $\left(\frac{\partial}{\partial z},\frac{\partial}{\partial z^*}\right)=\left(\frac{1}{2}\left(\frac{\partial}{\partial x}-\I \frac{\partial}{\partial y}\right),\frac{1}{2}\left(\frac{\partial}{\partial x}+\I \frac{\partial}{\partial y}\right)\right)$ we find from (\ref{theoremmellin})
\begin{eqnarray}
{\mit \Xi}(x+\I y) &=& {\hat{\mit \Omega}}\left\{1+\frac{1}{2}\left(x\frac{\partial}{\partial x}+y\frac{\partial}{\partial y}\right)+\frac{\I}{2}\left(y\frac{\partial}{\partial x}-x\frac{\partial}{\partial y}\right)\right\} \,{\rm ch}(x+\I y).\label{theoremmellinreal}
\end{eqnarray}
From this formula follows that ${\mit \Xi}(\I y)$ may be obtained by transformation of ${\rm ch}(\I y)=\cos(y)$ alone via
\begin{eqnarray}
{\mit \Xi}(\I y) &=& \hat{\mit \Omega}\left\{1+y\frac{\I}{2}\left(\frac{\partial}{\partial x}-\I \frac{\partial}{\partial y}\right) \right\}\,{\rm ch}\left(x+\I y\right) \;=\; \hat{\mit \Omega}\left(1+\I y\frac{\partial}{\partial z}\right)\,{\rm ch}(z) \nonumber\\ &=& \hat{\mit \Omega}\left(\frac{\partial}{\partial y}y\right)\cos(y).
\end{eqnarray}
On the right-hand side we have a certain redundance since in analytic functions the information which is contained in the values of the function on the imaginary axis is fully contained also in other parts of the function (here of ${\rm ch}(z)$).

The most simple transformation of ${\rm ch}(z)$ is by a delta function $\delta(u-u_0)$ as function ${\mit \Omega}(u)$ which stretches only the argument of the Hyperbolic Cosine function ${\rm ch}(z) \rightarrow {\rm ch}(u_0z)$. The next simple transformation is with a function
function ${\mit \Omega}(u)$ in form of a step function $\theta(u_0-u)$ which leads to the transformation ${\rm ch}(z) \rightarrow \frac{1}{z}\,{\rm sh}(u_0z)$. Our application of the second mean-value theorem reduced other cases under the suppositions of the theorem to this case, however, with parameter $u_0 =u_0(z)$ depending on complex variable $z$.

The great analogy between displacement operators (infinitesimal $-\I\frac{\partial}{\partial x}$ ) of the argument of a function and multiplication operator (infinitesimal $ x\frac{\partial}{\partial x}$) of the argument of a function with respect to the role of Fourier transformation and of Mellin transformation can be best seen from the following two relations
\begin{eqnarray}
\int_{-\infty}^{+\infty}dy\,f(y)g(x-y) &=& \tilde{f}\left(-\I \frac{\partial}{\partial x}\right)g(x),\quad \tilde{f}(t) \;\equiv \; \int_{-\infty}^{+\infty}dy\,f(y)\E^{-\I t y}, \nonumber\\ \int_{0}^{+\infty}du\,f(u)g(ux) &=& \hat{f}\left(\frac{\partial}{\partial x}x\right)g(x),\quad \hat{f}(s) \; \equiv \; \int_{0}^{+\infty}du\,f(u)\,u^{s-1}.\label{fourmellequ}
\end{eqnarray}
We remind you that Mellin and Fourier transform are related by substituting the integration variables $u=\E^y$ and the independent variables $s=-\I t$ and by the substitutions $f(\E^y) \leftrightarrow f(y)$ and $\hat{f}(-\I t) \leftrightarrow \tilde{f}(t)$ in (\ref{fourmellequ}).

Using the discussed Mellin transformation {\bf{Theorem 1}} can be reformulated as follows\\[1mm]
{\large{\bf Theorem $1^{\star}$:}}\\[1mm]
The mapping of the function ${\rm ch}(z)$ of the complex variable $z$ into the function ${\mit \Xi}(z)$ by an operator ${\hat{\mit \Omega}}\left(\frac{\partial}{\partial z}z\right)$
according to
\begin{eqnarray}
{\mit \Xi}(z) &=& {\hat{\mit \Omega}}\left(\frac{\partial}{\partial z}z\right) \,{\rm ch}(z),\quad
{\hat{\mit \Omega}}(s) \;\equiv\; \int_{0}^{+\infty}du\,{\mit \Omega}(u) u^{s-1},
\end{eqnarray}
where ${\hat{\mit \Omega}}(s)$ is the Mellin transformation of the function ${\mit \Omega}(u)$ which last possesses the properties given in {\bf{Theorem 1}} maps the function ${\rm ch}(z)$ with zeros only on the imaginary axis again into a function ${\mit \Xi}(z)$ with zeros only on the imaginary axis.\\
{\large{\bf Proof:}}\\[1mm]
It is proved as a reformulation of the {\bf{Theorem 1}} which is supposed here to be correctly proved.\\[2mm]
\indent
It was almost evident that the theorem may be formulated for more general functions ${\mit \Omega}(u)$ as supposed for the application of the second mean-value theorem that was made in Section 6. Under the suppositions of the theorem the integral on the left-hand side of (\ref{theoremmellin}) can be transformed by partial integration to (notation:
${\hat{\mit \Omega}^{(1)}}(s) \equiv \int_{0}^{+\infty}du\,{\mit \Omega}^{(1)}(u)u^{s-1}$)
\begin{eqnarray}
{\mit \Xi}(z) &=& -\frac{1}{z}\int_{0}^{+\infty}du\,{\mit \Omega}^{(1)}(u)\,{\rm sh}(uz) \;=\; -\frac{1}{z}\,{\hat{\mit \Omega}^{(1)}}\left(\frac{\partial}{\partial z}z\right)\,{\rm sh}(z).\label{theoremmellinsh}
\end{eqnarray}
The derivative ${\mit \Omega}^{(1)}(u)$ of the function ${\mit \Omega}(u)$ to the Riemann hypothesis although semi-definite (here negatively) and rapidly vanishing in infinity is not monotonic and possesses a minimum (see (\ref{defcapomega}) and Fig. 2.2).
In case of the (modified) Bessel functions we find by partial integration (e.g., \cite{erd})
\begin{eqnarray}
\nu!\left(\frac{2}{z}\right)^{\nu}\,{\rm I}_{\nu}\left(z\right) &=& \frac{\nu!} {\left(\nu-\frac{1}{2}\right)!\left(\frac{1}{2}\right)!}\int_{0}^{+1}du\,
\left(1-u^2\right)^{\nu-\frac{1}{2}}\,{\rm ch}\left(uz\right)\nonumber\\ &=& \frac{2\,\nu!}{\left(\nu-\frac{3}{2}\right)!\left(\frac{1}{2}\right)!}\,\frac{1}{z} \int_{0}^{+1}du\,\left(1-u^2\right)^{\nu-\frac{3}{2}}\,u\,{\rm sh}\left(uz\right),  \label{intreprbessel}
\end{eqnarray}
where the functions in the second transform $\left(1-u^2\right)^{\nu-\frac{3}{2}}\,u$ for $\nu > \frac{3}{2}$ are non-negative but not monotonic and possess a maximum for a certain value $u_{\rm max}$ within the interval $0 < u_{\rm max} <1$. The forms (\ref{theoremmellinsh}) for ${\mit \Xi}(z)$ and (\ref{intreprbessel}) suggest that there should be true a similar theorem to the integral in (\ref{basiccapxi2}) with substitution ${\rm ch}(uz) \rightarrow {\rm sh}(uz)$ and that monotonicity of the corresponding functions should not be the ultimate requirement for the zeros in such transforms on the imaginary axis.

Another consequence of the {\bf{Theorem 1}} follows from the non-negativity of the squared modulus of the function ${\mit \Xi}(z)$ resulting in the obvious inequality
\begin{eqnarray}
0 & \le & {\mit \Xi}(z)\left({\mit \Xi}(z)\right)^* \nonumber\\ &=& \frac{1}{2}\int_{0}^{+\infty}du_1\int_{0}^{+\infty}du_2\,{\mit \Omega}(u_1){\mit \Omega}(u_2) \nonumber\\ && \cdot \Big(\,{\rm ch}((u_1+u_2)x)\cos((u_1-u_2)y)+\,{\rm ch}((u_1-u_2)x)\cos((u_1+u_2)y)\Big),\label{ineq1}
\end{eqnarray}
which can be satisfied with the equality sign only on the imaginary axis $z=\I y$ for discrete values $y=y_k$ (the zeros of ${\mit \Xi}(z=x+\I y)$).
By transition from Cartesian coordinates $(u_1,u_2)$ to inertial-point coordinates $(u,\Delta u)$ according to
\begin{eqnarray}
&& u=\frac{u_1+u_2}{2},\quad \Delta u =u_2-u_1, \quad u_1=u-\frac{\Delta u}{2},\quad u_2= u+\frac{\Delta u}{2},\quad du_1 \wedge du_2 =du\wedge d(\Delta u),\qquad
\end{eqnarray}
equation (\ref{ineq1}) can be also written
\begin{eqnarray}
0 & \le & \frac{1}{2}\int_{0}^{+\infty}du\int_{-2u}^{+2u}d(\Delta u)\,{\mit
\Omega}\left(u-\frac{\Delta u}{2}\right){\mit \Omega}\left(u+\frac{\Delta u}{2}\right)
\nonumber\\ && \cdot
\Big(\,{\rm ch}\left(2u\,x\right)\cos\left(\Delta u\,y\right)+\,{\rm ch}\left(\Delta u\,x\right)\cos\left(2u\,y\right)\Big). \label{ineq2}
\end{eqnarray}
As already said the case of the equality sign in (\ref{ineq1}) or (\ref{ineq2}) can only be obtained for $x=0$ and then only for discrete values of $y=y_k$ by solution of this inequality with the specialization for $x=0$
\begin{eqnarray}
0 & = & {\mit \Xi}(\I y)\left({\mit \Xi}(\I y)\right)^* \nonumber\\ &=& \frac{1}{2}\int_{0}^{+\infty}du_1\int_{0}^{+\infty}du_2\,{\mit \Omega}(u_1){\mit \Omega}(u_2) \Big(\cos((u_1-u_2)y)+\cos((u_1+u_2)y)\Big) \nonumber\\ & = & \frac{1}{2}\int_{0}^{+\infty}du\int_{-2u}^{+2u}d(\Delta u)\,{\mit
\Omega}\left(u-\frac{\Delta u}{2}\right){\mit \Omega}\left(u+\frac{\Delta u}{2}\right)
\Big(\cos\left(\Delta u\,y\right)+\cos\left(2u\,y\right)\Big).\qquad \label{ineq3}
\end{eqnarray}
A short equivalent formulation of the inequality (\ref{ineq1}) and (\ref{ineq2}) together with (\ref{ineq3}) is the following\\
{\large{\bf Theorem 2:}}\\
If the function ${\mit \Omega}(u)$ satisfies the suppositions in {\bf{Theorem 1}} then with $y_{-k} = -y_{k}$
\begin{eqnarray}
{\mit \Xi}(x+\I y)\left({\mit \Xi}(x+ \I y)\right)^* \;\; \left\{\begin{array}{llll} > 0,& \Rightarrow & \left\{\begin{array}{cl} \forall x\neq 0,&  \\ \;\, x=0,& \forall y\neq y_k, \quad (k=\pm 1,\pm 2,\ldots), \end{array} \right. &
\\ =0, & \Rightarrow & \quad \;\;\, x=0, \quad \exists y=y_k, \quad (k=\pm 1,\pm 2,\ldots).\end{array} \right. \label{ineq4}
\end{eqnarray}
{\large{\bf Proof:}}\\[1mm]
As a consequence of proved {\bf{Theorem 1}} it is also proved. \\
The sufficient condition that this inequality is satisfied with the equality sign is that we first set $x=0$ in the expressions on the right-hand side of (\ref{ineq1}) and that we then determine the zeros $y=y_k$ of the obtained equation for ${\mit \Xi}(\I y)\left({\mit \Xi}(\I y)\right)^* =0$. In case of indefinite ${\mit \Omega}(u)$ there are possible in addition zeros on the $x$-axis.\\[1mm]
{\large{\bf Remark:}}\\[1mm]
Practically, (\ref{ineq1}) is an inequality for which it is difficult to prove in another way that it can be satisfied with the equality sign only for $x=0$. Proved in another way with specialization (\ref{defcapomega}) for ${\mit \Omega}(u)$ it would be an independent proof of the Riemann hypothesis.

\setcounter{chapter}{9}
\setcounter{equation}{0}
\section*{Conclusion}

We proved in this article the Riemann hypothesis embedded into a more general theorem for a class of functions ${\mit \Xi}(z)$ with a representation of the form (\ref{basiccapxi0}) for real-valued functions ${\mit \Omega}(u)$ which are positive semi-definite and non-increasing in the interval $0 \le u < +\infty$ and which are vanishing in infinity more rapidly than any exponential function $\exp\left(-\lambda u\right)$ with $\lambda > 0$. The special Xi function ${\mit \Xi}(z)$ to the function ${\mit \Omega}(u)$ given in (\ref{defcapomega}) which is essentially the xi function $\xi(s)$ equivalent to the Riemann zeta function $\zeta(s)$ concerning the hypothesis belongs to the described class of functions.

Modified Bessel functions of imaginary argument 'normalized' to entire functions $(\frac{2}{\I z})^{\nu}\,{\rm J}_{\nu}(\I z) =(\frac{2}{z})^{\nu}\,{\rm J}_{\nu}(z)$ for $\nu \ge \frac{1}{2}$ belong also to this class of functions with a representation of the form (\ref{basiccapxi0}) with ${\mit \Omega}(u)$ which satisfy the mentioned conditions and in this last case it is well known and proved in independent way that their zeros lie only on the imaginary axis corresponding to the critical line in the Riemann hypothesis. Knowing this property of the modified Bessel functions we looked from beginning for whole classes of functions including the Riemann zeta function which satisfy analogous conditions as expressed in the Riemann hypothesis. The details of the approach to Bessel functions and also to certain classes of almost-periodic functions we prepare for another work.

The numerical search for zeros of the Riemann zeta function $\zeta(s),\,s=\sigma + \I t$ in the critical strip, in particular, off the critical line may come now to an end by the proof of the Riemann hypothesis since its main purpose was, in our opinion, to find a counter-example to the Riemann hypothesis and thus to disprove it. We did not pay attention in this article to methods of numerical calculation of the zeros with (ultra-)high precision and for very high values of the imaginary part. However, the proof if correct may deliver some calculators now from their pain to calculate more and more zeros of the Riemann zeta function.

We think that some approaches in this article may possess importance also for other problems. First of all this is the operational approach of the transition from real and imaginary part of a function on the real or imaginary axis to an analytic function in the whole complex plane. In principle, this is possible using the Cauchy-Riemann equations but the operational approach integrates this to two integer instead of differential equations. We think that this is possible also in curved coordinates and is in particular effective starting from curves of constant real or imaginary part of one of these functions on a curve.

One of the fascinations of prime number theory is the relation of the apparently chaotic distribution function of prime numbers $\pi(x)$ on the real axis $x\ge 0$ to a fully well-ordered analytic function, the Riemann zeta function $\zeta(s)$, at least, in its representation in sum form as a special Dirichlet series and thus providing the relations between multiplicative and additive representations of arithmetic functions.

{\setcounter{section}{20}}
\section*{Acknowledgement}
I am very grateful to my son Arne W\"{u}nsche for different help. He helped me very much with the arrangement of my computer and if there were problems with it and also with the installation and use of programmes and by often simple questions he gave me new suggestions.


{\setcounter{chapter}{21}
\setcounter{equation}{0}
\aeqn
\section*{Appendix A.\\ Transformation of the Xi function\\}

In this Appendix we transform the function $\xi(s)$ defined in (\ref{xirep1}) by means of the zeta function $\zeta(s)$ from the form taken from (\ref{zetarep2}) to the form (\ref{xirep2}) using the Poisson summation formula. The Poisson summation formula is the transformation of a sum over a lattice into a sum over the reciprocal lattice. More generally, in one-dimensional case the decomposition of a special periodic function $F(q)=F(q+a)$ with period $a$ defined by the following series over functions $f(q+na)$
\begin{eqnarray}
F(q) \;\equiv\; \sum_{n=-\infty}^{+\infty}f(q+na) \;=\; F(q+a),
\end{eqnarray}
can be transformed into the reciprocal lattice providing a Fourier series as follows. For this purpose we expand $F(q)$ in a Fourier series with Fourier coefficients $F_m$ and make then obvious transformations ($q'+na \rightarrow q'',\,\dis{\sum_{n=-\infty}^{+\infty}\int_{na}^{(n+1)a}dq'' \rightarrow \int_{-\infty}^{+\infty}dq''} $ and changing the order of summation and integration) according to
\begin{eqnarray}
F(q) &=& \sum_{m=-\infty}^{+\infty}\underbrace{\left(\frac{1}{a}\int_{0}^{a}dq'\sum_{n=-\infty}^{+\infty}f(q'+na)\exp\left(-\I m\frac{2\pi q'}{a}\right)\right)}_{=\,F_m}\exp\left(\I m \frac{2\pi q}{a}\right) \nonumber\\ &=&
\frac{1}{a}\sum_{m=-\infty}^{+\infty}\left(\sum_{n=-\infty}^{+\infty}\int_{na}^{(n+1)a}dq''\,f(q'')\exp\left(-\I m\frac{2\pi q''}{a}\right)\right) \exp\left(\I m\frac{2\pi q}{a}\right) \nonumber\\ &=& \frac{1}{a}\sum_{m=-\infty}^{+\infty} \tilde{f}\left(m\frac{2\pi}{a} \right)\exp\left(\I m\frac{2\pi q}{a}\right),\label{poisson1}
\end{eqnarray}
where the coefficients $F_{m}$ of the decomposition of $F(q)$ are given by the Fourier transform $\tilde{f}(k)$ of the function $f(q)$ defined in the following way
\begin{eqnarray}
\tilde{f}(k) \;=\; \int_{-\infty}^{+\infty}dx f(q)\exp\left(-\I kq\right),\quad \Rightarrow \quad
\tilde{f}\left(m\frac{2\pi}{a}\right) \;=\; \int_{-\infty}^{+\infty}dq f(q)\exp\left(-\I m \frac{2\pi q}{a}\right).\label{fourcoeff}
\end{eqnarray}
Using the period $b= \frac{2\pi}{a}$ of the reciprocal lattice relation on the right-hand side of (\ref{poisson1}) it may be written in the forms
\begin{eqnarray}
\sum_{n=-\infty}^{+\infty}f(q+na) &=& \sum_{n=-\infty}^{+\infty}f\left(q+n\frac{2\pi}{b}\right)\;=\; \frac{b}{2\pi}\sum_{m=-\infty}^{+\infty} \tilde{f}\left(mb\right)\exp\left(\I mb q\right) \nonumber\\ &=&  \frac{1}{a}\sum_{m=-\infty}^{+\infty} \tilde{f}\left(m\frac{2\pi}{a} \right)\exp\left(\I m\frac{2\pi}{a}q\right), \qquad ab = 2\pi.\qquad \label{poisson3}
\end{eqnarray}
In the special case $q=0$ one obtains from (\ref{poisson1}) the well-known basic form of the Poisson summation formula
\begin{eqnarray}
F(0) \;\equiv \; \sum_{n=-\infty}^{+\infty}f(na) \;=\;
\frac{b}{2\pi}\sum_{m=-\infty}^{+\infty}\tilde{f}\left(mb\right), \qquad ab=2\pi.\label{poisson2}
\end{eqnarray}

Formula (\ref{poisson2}) applied to the sum $\dis{\sum_{n=-\infty}^{+\infty}\exp\left(-\pi n^2q^2\right)}$ corresponding to $f(q)=\exp\left(-\pi q^2\right)$ with Fourier transform $\tilde{f}(k)=\exp\left(-\frac{k^2}{4\pi}\right)$ provides a relation which can be written in the following symmetric form (we need it in the following only for $q\ge 0$)
\begin{eqnarray}
{\mit \Psi}(q) &\equiv & \sqrt{\left|q\right|}\left\{\frac{1}{2}+\sum_{n=1}^\infty
\exp\left(-\pi n^2 q^2\right)\right\} \;=\;\frac{\sqrt{\left|q\right|}}{2} \sum_{n=-\infty}^{+\infty} \exp\left(-\pi n^2 q^2\right) \nonumber\\ &=& \frac{1}{2\sqrt{\left|q\right|}} \sum_{m=-\infty}^{+\infty}
\exp\left(-\pi \frac{m^2}{q^2}\right) \;=\; \frac{1}{\sqrt{\left|q\right|}} \left\{\frac{1}{2}+\sum_{m=1}^\infty \exp\left(-\pi \frac{m^2}{q^2}\right)\right\} \;\equiv\; {\mit \Psi}\left(\frac{1}{q}\right).
\label{defcappsi}
\end{eqnarray}
This is essentially a transformation of the Theta function $\vartheta_3\left(u,q\right)$ in special case $\frac{2{\mit \Psi}(q)}{\scr{\sqrt{\left|q\right|}}} \equiv \vartheta_3\left(0,\E^{-\pi q^2}\right)$. We now apply this to a transformation of the function $\xi(s)$.

From (\ref{xirep2}) and (\ref{zetarep2}) follows
\begin{eqnarray}
\xi(s) &=& s(s-1)\int_{0}^{+\infty}dq\,q^{s-1} \sum_{n=1}^\infty\exp\left(-\pi n^2q^2\right) \nonumber\\ &=& s(s-1)\left\{\int_{0}^{1}dq\,q^{s-1} \sum_{n=1}^\infty\exp\left(-\pi
n^2q^2\right)+\int_{1}^{+\infty}dq\,q^{s-1} \sum_{n=1}^\infty\exp\left(-\pi n^2q^2\right) \right\}.\label{xi1}
\end{eqnarray}
The second term in braces is convergent for arbitrary $q$ due to the rapid vanishing of the summands of the sum for $q\rightarrow \infty$. To the first term in braces we apply the Poisson summation formula (\ref{poisson2}) and obtain from the special result
(\ref{defcappsi})
\begin{eqnarray}
\int_{0}^{1}dq\,q^{s-1} \sum_{n=1}^\infty\exp\left(-\pi n^2q^2\right) &=& \int_{0}^{1}dq\,q^{s-1} \left\{\frac{1}{2}\left(\frac{1}{q}-1\right)+ \frac{1}{q}\sum_{m=1}^{\infty}
\exp\left(-\pi\frac{m^2}{q^2}\right)\right\}\nonumber\\ &=& \frac{1}{2s(s-1)}+\int_{1}^{+\infty}dq'\, q'^{-s}\sum_{m=1}^\infty \exp\left(-\pi m^2 q'^2\right),
\end{eqnarray}
with the substitution $q=\frac{1}{q'}$ of the integration variable made in last line. Thus from (\ref{xi1}) we find
\begin{eqnarray}
\xi(s) &=&
\frac{1}{2}-s(1-s)\int_{1}^{+\infty}dq\,\frac{q^{s}+q^{1-s}}{y} \sum_{n=1}^\infty\exp\left(-\pi n^2q^2\right).
\end{eqnarray}
With the substitution of the integration variable
\begin{eqnarray}
q=\E^u,\quad \left(q\ge 0,\; \Leftrightarrow \; -\infty < u <+\infty\right),\label{yex}
\end{eqnarray}
and with displacement of the complex variable $s$ to $z \equiv s+\frac{1}{2}$ and introduction of ${\mit \Xi}\left(z\right)$ instead of $\xi(s)$ this leads to the representation
\begin{eqnarray}
{\mit \Xi}\left(z\right)&\equiv & \xi\left(\frac{1}{2}+z\right)\;=\; \frac{1}{2}
-2\left(\frac{1}{4}-z^2\right)\int_{0}^{+\infty}du\,{\rm ch}\left(u z\right)\,\E^{\frac{u}{2}} \sum_{n=1}^\infty \exp\left(-\pi n^2 \E^{2u}\right),\label{capxirepa3}
\end{eqnarray}
given in (\ref{capxirep2}). In Appendix B we transform this representation by means of partial integration to a form which due to symmetries is particularly appropriate for the further
considerations about the Riemann zeta function.

Using the substitution (\ref{yex}) we define a function ${\mit \Phi}(u)$ by means of the function ${\mit \Psi}(y)$ in (\ref{defcappsi}) as follows
\begin{eqnarray}
{\mit \Phi}\left(u\right) & \equiv & {\mit \Psi}\left(\E^{u}\right),\quad \Rightarrow \quad {\mit
\Phi}\left(0\right) \;=\; {\mit \Psi}\left(1\right) \;=\; 0.543217\label{defcapphi},
\end{eqnarray}
and explicitly due to Poisson summation formula
\begin{eqnarray}
{\mit \Phi}\left(u\right) & \equiv & \E^{\frac{u}{2}}\left\{\frac{1}{2}+\sum_{n=1}^\infty \exp\left(-\pi n^2 \E^{2u}\right)\right\} \;=\; \frac{1}{2}\E^{\frac{u}{2}}\sum_{n=-\infty}^{+\infty} \exp\left(-\pi
n^2 \E^{2u}\right) \nonumber\\ &=& \frac{1}{2}\E^{-\frac{u}{2}}\sum_{n=-\infty}^{+\infty}
\exp\left(-\pi n^2 \E^{-2u}\right) \;=\; \E^{-\frac{u}{2}}\left\{\frac{1}{2}+\sum_{n=1}^\infty \exp\left(-\pi n^2 \E^{-2u}\right)\right\}.
\end{eqnarray}
From ${\mit \Psi}(y) = {\mit \Psi}\!\left(\frac{1}{y}\right)$ according to (\ref{defcappsi}) follows that ${\mit \Phi}\left(u\right)$ is a symmetric function
\begin{eqnarray}
{\mit \Phi}\left(u\right) = {\mit \Phi}\left(-u\right) = {\mit \Phi}\left(\left|u\right|\right).
\end{eqnarray}
Therefore, all even derivatives of ${\mit \Phi}(u)$ are also symmetric functions, whereas all odd derivatives of ${\mit \Phi}(u)$ are antisymmetric functions (we denote these derivatives by ${\mit
\Phi}^{\left(n\right)}(u)\equiv \frac{\partial^n {\mit \Phi}}{\partial u}(u)$)
\begin{eqnarray}
{\mit \Phi}^{\left(2m\right)}\left(u\right) &=& +{\mit \Phi}^{\left(2m\right)}\left(-u\right)\;=\; {\mit \Phi}^{\left(2m\right)}\left(\left|u\right|\right),\nonumber\\ {\mit \Phi}^{\left(2m+1\right)}\left(u\right) &=& -{\mit \Phi}^{\left(2m+1\right)}\left(-u\right),\quad \Rightarrow \quad {\mit \Phi}^{\left(2m+1\right)}\left(0\right)\;=\;0, \quad \left(m=0,1,2,\ldots\right).\label{symmcapphi}
\end{eqnarray}
Explicitly, one obtains for the first two derivatives
\begin{eqnarray}
{\mit \Phi}^{\left(1\right)}\left(u\right) &=& \frac{1}{2}\E^{\frac{u}{2}}\left\{\frac{1}{2}+ \sum_{n=1}^\infty \left(1-4\pi n^2\E^{2u}\right) \exp\left(-\pi n^2 \E^{2u}\right)\right\},\nonumber\\ {\mit \Phi}^{\left(2\right)}\left(u\right) &=&
\frac{1}{4}\E^{\frac{u}{2}}\left\{\frac{1}{2}+ \sum_{n=1}^\infty \left(1-6\left(4\pi n^2\E^{2u}\right)+\left(4\pi n^2\E^{2u}\right)^2\right) \exp\left(-\pi n^2 \E^{2u}\right)\right\}.
\label{derivcapphi}
\end{eqnarray}
As a subsidiary result we obtain from vanishing of the odd derivatives of ${\mit \Phi}(u)$ at $u=0$ that means from ${\mit \Phi}^{\left(2m+1\right)}\left(0\right)=0,\,(m=0,1,2,\ldots)$ an infinite sequence of special sum evaluations from which the first three are
\begin{eqnarray}
&& \hspace{-5mm}\sum_{n=1}^\infty \left(4\,\pi n^2-1\right)\exp\left(-\pi n^2\right) \;=\; \frac{1}{2},\nonumber\\ && \hspace{-5mm} \sum_{n=1}^\infty \left(\left(4\,\pi n^2\right)^3-15\left(4\,\pi n^2\right)^2 +31\left(4\,\pi n^2\right)-1\right)\exp\left(-\pi n^2\right) \;=\; \frac{1}{2}.
\label{derivcapphizero}
\end{eqnarray}
We checked relations (\ref{derivcapphizero}) numerically by computer up to a sufficiently high precision. We also could not find (\ref{gentheta}) among the known transformations of theta functions. The interesting feature of these sum evaluations is that herein power functions as well as exponential functions containing the transcendental number $\pi$ in the exponent are involved in a way which finally leads to a rational number that should also be attractive for recreation mathematics. In contrast, in the well-known series for the trigonometric functions one obtains for certain rational multiples of $\pi$ as argument also rational numbers but one has involved there only power functions with rational coefficients that means rational functions although an infinite number of them.

Using the function ${\mit \Phi}(u)$ the function ${\mit \Xi}(z)$ in (\ref{capxirepa3}) can be represented as
\begin{eqnarray}
{\mit \Xi}\left(z\right) &=& \frac{1}{2} +2\left(z^2-\frac{1}{4}\right)\int_{0}^{+\infty}du\,\left({\mit \Phi}\left(u\right)-\frac{1}{2}\,\E^{\frac{u}{2}}\right)\,{\rm ch}\left(uz\right) \nonumber\\ &=& \frac{1}{2} +2\int_{0}^{+\infty}du\,\left({\mit \Phi}\left(u\right)-\frac{1}{2}\,\E^{\frac{u}{2}}\right) \left\{\left(\frac{\partial^2}{\partial
u^2}-\frac{1}{4}\right)\,{\rm ch}\left(uz\right)\right\}.\label{capxirepa4}
\end{eqnarray}
From this we obtain by partial integration
\begin{eqnarray}
{\mit \Xi}\left(z\right) &=& 2\int_{0}^{+\infty}du\,\left\{\left(\frac{\partial^2}{\partial
u^2}-\frac{1}{4}\right)\left({\mit \Phi}\left(u\right)-\frac{1}{2}\,\E^{\frac{u}{2}}\right)\right\}\,{\rm ch}\left(uz\right),\label{capxirepa5}
\end{eqnarray}
where the contribution from the lower integration limit at $u=0$ has exactly canceled the constant term $\frac{1}{2}$ on the right-hand side of (\ref{capxirepa4}) and the contributions from the upper limit $x\rightarrow +\infty$ is vanishing. Using (\ref{derivcapphi}) we find with abbreviation ${\mit \Omega}(u)$ according to
\begin{eqnarray}
{\mit \Omega}(u) &\equiv & 2\left(\frac{\partial^2}{\partial u^2}-\frac{1}{4}\right)\left({\mit
\Phi}\left(u\right)-\frac{1}{2}\,\E^{\frac{u}{2}}\right),
\end{eqnarray}
the following basic structural form of the Xi function
\begin{eqnarray}
{\mit \Xi}\left(z\right) &=& \int_{0}^{+\infty}du\,{\mit \Omega}(u)\,{\rm ch}\left(uz\right),\label{capxirepa6}
\end{eqnarray}
with the following explicit representation of ${\mit \Omega}(u)$
\begin{eqnarray}
{\mit \Omega}(u) &=& 2{\mit \Phi}^{\left(2\right)}\left(u\right)-\frac{1}{2}\,{\mit
\Phi}\left(u\right) \nonumber\\ &=& 4 \E^{\frac{u}{2}}\sum_{n=1}^\infty \pi n^2\E^{2u}\left(2\pi n^2 \E^{2u}-3\right)\exp\left(-\pi n^2\E^{2u}\right) \;>\;0.\label{defcapomega1}
\end{eqnarray}
Since according to (\ref{symmcapphi}) the even derivatives of ${\mit \Phi}(u)$ are symmetric functions it follows from relation (\ref{defcapomega1}) that ${\mit \Omega}(u)$ is also a symmetric function and (\ref{omegasym}) holds. This is not immediately seen from the explicit
representation (\ref{defcapomega1}). Furthermore, ${\mit \Omega}(u)$ is positively definite for $-\infty < u < +\infty $ since the factor $\left(2\pi n^2 \E^{2u}-3\right)$ in
(\ref{defcapomega1}) is positive for $n\ge 1$ and $u\ge 0$ and all other factors too. It goes rapidly to zero for $u\rightarrow \pm\infty$, more rapidly than any exponential function $\exp\left(-\gamma|x|^\rho\right)$ with arbitrary $\gamma>0$ and arbitrary $\rho >0$ due to factors $\exp\left(-\pi n^2 \E^{2u}\right)$ in the sum terms in (\ref{defcapomega1}). For
the first derivative of ${\mit \Omega}(u)$ we find
\begin{eqnarray}
{\mit \Omega}^{\left(1\right)}(u) &=& 2{\mit \Phi}^{\left(3\right)}\left(u\right)-\frac{1}{2}\,{\mit \Phi}^{\left(1\right)}\left(u\right) \nonumber\\ &=& -2 \E^{\frac{u}{2}}\sum_{n=1}^\infty \pi n^2\E^{2u}\left(8\left(\pi
n^2\E^{2u}\right)^2-30\pi n^2\E^{2u}+15\right)\exp\left(-\pi n^2 \E^{2u}\right) \nonumber\\ &=& -{\mit \Omega}^{\left(1\right)}(-u), \quad \Rightarrow \quad {\mit \Omega}^{\left(1\right)}(0)=0,\quad {\mit \Omega}^{\left(1\right)}(u) \;<\; 0,\quad \left(u > 0\right).
\end{eqnarray}
It is vanishing for $u=0$ due to its antisymmetry and negatively definite for $u > 0$ as the negative sign of ${\mit \Omega}^{\left(2\right)}(0)$ together with considerations of the sum
for $u>0$ show (i.e., the polynomial $f(N)\equiv 8\pi^2 N^2-30\pi N+15 \ge 0 $ for $N\ge \frac{15+\sqrt{105}}{8\pi} \approx 1.00454 $ and negativity is already obtained taking the first two sum terms to $n=1$ and $n=2$ alone). Thus ${\mit \Omega}(u)$ is monotonically  decreasing for $u\ge 0$. A few approximate numerical values of parameters for the function ${\mit \Omega}(u)$ are
\begin{eqnarray}
&& {\mit \Omega}_0 \;\equiv \; \int_{0}^{+\infty}du\,{\mit \Omega}(u) \;=\; -\int_{0}^{+\infty}du\,{\mit \Omega}^{(1)}(u)u \;=\; 0.497121,\nonumber\\ && {\mit \Omega}(0) \;=\; - \int_{0}^{+\infty}dx\,{\mit \Omega}^{(1)}(u) \;=\; 1.78679,\quad {\mit \Omega}^{\left(1\right)}(0) \;=\; - \int_{0}^{+\infty}du\,{\mit \Omega}^{(2)}(u) \;=\; 0.
\label{numcapomega}
\end{eqnarray}

In next Appendix we consider the transition from analytic functions given on the real or imaginary axis to the whole complex plane.

}

{\setcounter{chapter}{15}
\beqn
\setcounter{equation}{0}
\section*{Appendix B.\\ Transition from analytic functions on real or imaginary axis to whole complex plane}

The operator $\frac{\partial}{\partial x}$ is the infinitesimal displace operator and $\exp\left(-x_0\frac{\partial}{\partial x}\right)$ the finite displacement operator for the displacement of the argument of a function $f(x) \rightarrow f(x-x_0)$. In complex analysis the
real variable $x$ can be displaced with sense for an analytic function to the complex variable $z=x+\I y$ in the whole complex plane by
\begin{eqnarray}
\exp\left(\I y \frac{\partial}{\partial x}\right)x \exp\left(-\I y \frac{\partial}{\partial x}\right) &=& x+\frac{\I y}{1!}\left[\frac{\partial}{\partial x},x \right] + \frac{(\I y)^2}{2!}\left[\frac{\partial}{\partial x}\left[\frac{\partial}{\partial x},x\right]\right]+\ldots \;=\; x+\I y, \qquad\label{xdisply0}
\end{eqnarray}
where $\left[A,B\right] \equiv AB-BA $ denotes the commutator of two operators $A$ and $B$, in particular $\left[\frac{\partial}{\partial x},x\right] =1$ and (\ref{xdisply0}) may be written in the form
\begin{eqnarray}
\exp\left(\I y \frac{\partial}{\partial x}\right)x &=& (x+\I y)\exp\left(\I y \frac{\partial}{\partial x}\right).\label{xdisply}
\end{eqnarray}
Analogously, the transition from the variable $y$ on the imaginary axis $\I y$ to the variable $z=\I(y-\I x)=x+\I y $ in the whole complex plane may be written as
\begin{eqnarray}
\exp\left(-\I x \frac{\partial}{\partial y}\right)\I y &=& (x+\I y)\exp\left(-\I x \frac{\partial}{\partial y}\right),\qquad \; \label{ydisplx}
\end{eqnarray}
In the following we consider only the case (\ref{xdisply}) since the case (\ref{ydisplx}) is
completely analogous with simple substitutions.

We wrote the equations (\ref{xdisply0}), (\ref{xdisply}) and (\ref{ydisplx}) in a form which we call operational form and meaning that they may be applied to further functions on the left-hand and correspondingly right-hand side\footnote{Non-operational form would be if we write, for example, $\exp\left(\I y \frac{\partial}{\partial x}\right)x = x+\I y$  instead of (\ref{xdisply}) which is correct but cannot be applied to further functions $f(x) \neq const\cdot 1 $, for example to $f(x)=x$.)}. It is now easy to see that an analytic function $w(z)=w(x+\I y); \frac{\partial}{\partial z^*}w(z) =0$ can be generated from $w(x)$ on the $x$-axis in operational form by
\begin{eqnarray}
\exp\left(\I y \frac{\partial}{\partial x}\right)w(x) &=& w(x+\I y)\exp\left(\I y \frac{\partial}{\partial x}\right),\label{fxdisplfxy}
\end{eqnarray}
and analogously from $w(\I y)$ on the imaginary axis by
\begin{eqnarray}
\exp\left(-\I x \frac{\partial}{\partial x}\right)w(\I y) &=& w(x+\I y)\exp\left(-\I x \frac{\partial}{\partial y}\right).\label{fydisplfxy}
\end{eqnarray}

Writing the function $w(z)$ with real part $u(x,y)$ and imaginary part $v(x,y)$ in the form
\begin{eqnarray}
&& w(x+\I y)= u(x,y)+\I v(x,y),\quad \left(w(x+\I y)\right)^* = u(x,y)-\I v(x,y), \quad \Rightarrow \nonumber\\&& u(x,y) \;=\; \frac{1}{2}\Big(w(x+\I y)+ \left(w(x+\I y)\right)^*\Big),\quad v(x,y) \;=\; -\frac{\I}{2}\Big(w(x+\I y)-\left(w(x+\I y)\right)^*\Big),\qquad
\end{eqnarray}
we find from (\ref{fxdisplfxy})
\begin{eqnarray}
\exp\left(\I y \frac{\partial}{\partial x}\right)w(x) &=& \exp\left(\I y \frac{\partial}{\partial x}\right)\Big(u(x,0)+\I v(x,0)\Big) \nonumber\\ &=& \Big(u(x+\I y,0)+ \I v(x+\I y,0)\Big)\exp\left(\I y \frac{\partial}{\partial x}\right) \nonumber\\ &=& \Big(u(x,y)+ \I v(x,y)\Big)\exp\left(\I y \frac{\partial}{\partial x}\right),\label{fxdisplfxyplus}
\end{eqnarray}
and correspondingly
\begin{eqnarray}
\left(\exp\left(\I y \frac{\partial}{\partial x}\right)w(x)\right)^* &=& \exp\left(-\I y \frac{\partial}{\partial x}\right)\Big(u(x,0)-\I v(x,0)\Big) \nonumber\\ &=& \Big(u(x-\I y,0)- \I v(x-\I y,0)\Big)\exp\left(-\I y \frac{\partial}{\partial x}\right) \nonumber\\ &=& \Big(u(x,y)- \I v(x,y)\Big)\exp\left(-\I y \frac{\partial}{\partial x}\right).\label{fxdisplfxyminus}
\end{eqnarray}
From (\ref{fxdisplfxyplus}) and (\ref{fxdisplfxyminus}) follows forming the sum and the difference
\begin{eqnarray}
\cos\left(y \frac{\partial}{\partial x}\right)u(x,0) -\sin\left(y \frac{\partial}{\partial x}\right)v(x,0) &=& u(x,y)\cos\left(y \frac{\partial}{\partial x}\right) -v(x,y)\sin\left(y \frac{\partial}{\partial x}\right), \nonumber\\ \sin\left(y \frac{\partial}{\partial x}\right)u(x,0) +\cos\left(y \frac{\partial}{\partial x}\right)v(x,0) &=& u(x,y)\sin\left(y \frac{\partial}{\partial x}\right) +v(x,y)\cos\left(y \frac{\partial}{\partial x}\right).\qquad
\label{fxdisplfxyop}
\end{eqnarray}
These are yet operational identities which can be applied to arbitrary functions $f(x)$. Applied to the function $f(x)=1$ follows
\begin{eqnarray}
\cos\left(y \frac{\partial}{\partial x}\right)u(x,0) -\sin\left(y \frac{\partial}{\partial x}\right)v(x,0) &=& u(x,y), \nonumber\\ \sin\left(y \frac{\partial}{\partial x}\right)u(x,0) +\cos\left(y \frac{\partial}{\partial x}\right)v(x,0) &=& v(x,y).
\label{fxdisplfxynonop}
\end{eqnarray}

In full analogy we may derive the continuation of an analytic function from the imaginary axes $z=\I y$ to the whole complex plane $z=x+\I y$ in operational form
\begin{eqnarray}
\cos\left(x\frac{\partial}{\partial y}\right)u(0,y) +\sin\left(x\frac{\partial}{\partial y}\right)v(0,y) &=& u(x,y)\cos\left(x\frac{\partial}{\partial y}\right) +v(x,y)\sin\left(x\frac{\partial}{\partial y}\right), \nonumber\\ -\sin\left(x\frac{\partial}{\partial y}\right)u(0,y) +\cos\left(x\frac{\partial}{\partial y}\right)v(0,y) &=& -u(x,y)\sin\left(x\frac{\partial}{\partial y}\right) +v(x,y)\cos\left(x\frac{\partial}{\partial y}\right),\qquad \label{fydisplfxyop}
\end{eqnarray}
and this applied to the function $f(y)=1$
\begin{eqnarray}
\cos\left(x\frac{\partial}{\partial y}\right)u(0,y) +\sin\left(x\frac{\partial}{\partial y}\right)v(0,y) &=& u(x,y), \nonumber\\ -\sin\left(x\frac{\partial}{\partial y}\right)u(0,y) +\cos\left(x\frac{\partial}{\partial y}\right)v(0,y) &=& v(x,y). \label{fydisplfxynonop}
\end{eqnarray}
It is easy to check that both (\ref{fxdisplfxynonop}) and (\ref{fydisplfxynonop}) satisfy the Cauchy-Riemann equations
\begin{eqnarray}
\frac{\partial u}{\partial x}(x,y) \;=\; \frac{\partial v}{\partial y}(x,y),\quad \frac{\partial u}{\partial y}(x,y) \;=\; -\frac{\partial v}{\partial x}(x,y),\label{creq}
\end{eqnarray}
and it is even possible to derive these relations from these equations by Taylor series expansions of $u(x,y)$ and $v(x,y)$ in powers of $y$ or $x$ in dependence from which axis we make the continuation to the whole complex plane. For example, in expansion in powers of $y$ we obtain
using (\ref{creq}) (and the resulting equations $\left(\frac{\partial^2}{\partial x^2}+\frac{\partial^2}{\partial y^2}\right)u(x,y) =\left(\frac{\partial^2}{\partial x^2}+\frac{\partial^2}{\partial y^2}\right)v(x,y) = 0$ from them)
\begin{eqnarray}
w(x+\I y) &=& \sum_{m=0}^{\infty}\frac{y^{2m}}{(2m)!}\left(\frac{\partial^{2m}}{\partial y^{2m}}u(x,y)\right)_{y=0} + \sum_{m=0}^{\infty}\frac{y^{2m+1}}{(2m+1)!}\left(\frac{\partial^{2m+1}}{\partial y^{2m+1}}u(x,y)\right)_{y=0} \nonumber\\ && +\I\left\{\sum_{m=0}^{\infty}\frac{y^{2m}}{(2m)!}\left(\frac{\partial^{2m}}{\partial y^{2m}}v(x,y)\right)_{y=0} + \sum_{m=0}^{\infty}\frac{y^{2m+1}}{(2m+1)!}\left(\frac{\partial^{2m+1}}{\partial y^{2m+1}}v(x,y)\right)_{y=0} \right\},\qquad\quad
\end{eqnarray}
that can be written in compact form
\begin{eqnarray}
w(x+\I y) &=& \cos\left(y\frac{\partial}{\partial x}\right)u(x,0)-\sin\left(y\frac{\partial}{\partial x}\right)v(x,0) \nonumber\\ && +\I\left\{\cos\left(y\frac{\partial}{\partial x}\right)v(x,0)+\sin\left(y\frac{\partial}{\partial x}\right)u(x,0)\right\}\;=\; u(x,y) +\I v(x,y),\label{fxdisplfxynonop2}
\end{eqnarray}
and is equivalent to (\ref{fxdisplfxynonop}). Analogously by expansion in powers of $x$ as intermediate step we obtain
\begin{eqnarray}
w(x+\I y) &=& \cos\left(x\frac{\partial}{\partial y}\right)u(0,y)+\sin\left(x\frac{\partial}{\partial y}\right)v(0,y) \nonumber\\ && +\I\left\{\cos\left(x\frac{\partial}{\partial y}\right)v(0,y)-\sin\left(x\frac{\partial}{\partial y}\right)u(0,y)\right\} \;=\; u(x,y)+\I v(x,y),\label{fydisplfxyop2}
\end{eqnarray}
that is equivalent to (\ref{fydisplfxyop}).
Therefore, relations (\ref{fxdisplfxynonop2}) and (\ref{fydisplfxyop2}) represent some integral forms of the Cauchy-Riemann equations.

In cases if one of the functions $u(x,0)$ or $v(x,0)$ in (\ref{fxdisplfxynonop}) or $u(0,y)$ or $v(0,y)$ in (\ref{fydisplfxynonop}) is vanishing these formulae simplify and the case $v(0,y) = 0$ is applied in Sections 4--6.
We did not find up to now such representations in textbooks to complex analysis but it seems possible that they are somewhere.

}

\small{


\end{document}